\documentstyle[12pt]{article}
\def\Bbb R{{\rm \bf R}}
\def\proclaim#1{\vskip2mm{\bf #1}\em}
\def\endproclaim{\em \vskip2mm}
\def\tag#1{\eqno(#1)}
\def\gathered{\begin{array}{c}}
\def\endgathered{\end{array}}
\def\text{\mbox}

\begin{document}

\title {The enclosure method for inverse obstacle scattering problems with dynamical data
over a finite time interval II.  Obstacles with a dissipative boundary or
finite refractive index and back-scattering data}
\author{Masaru IKEHATA\footnote{
Department of Mathematics,
Graduate School of Engineering,
Gunma University, Kiryu 376-8515, JAPAN}}
\maketitle

\begin{abstract}
In this paper a wave is generated by an initial data whose support
is localized at the outside of unknown obstacles and observed in a {\it limited time}
on a known closed surface or the same position as the support of
the initial data. The observed data in the latter process are
nothing but the {\it back-scattering data}. Two types of obstacles
are considered.  One is obstacles with a {\it dissipative boundary
condition} which is a generalization of the {\it sound-hard
obstacles}; another is obstacles with a {\it finite refractive
index}, so-called, transparent obstacles. For each type of
obstacles two formulae which yield explicitly the distance
from the support of the initial data to unknown obstacles are
given.

\noindent
AMS: 35R30

\end{abstract}

\section{Introduction}

In \cite{IE0} the author has introduced a simple method for the
{\it reconstruction issue} of some class of inverse obstacle
scattering problems. The observation data used therein are given
by a wave field observed in limited time on a known closed surface
surrounding unknown obstacles. The wave is generated by the
initial data with compact support outside the surface. It is shown
that the method yields the distance from a given point outside the
surface to obstacles provided the obstacles are {\it sound-hard}
or {\it penetrable}. It is a kind of an {\it enclosure method}
since the obstacles can be considered as being contained in an
open ball centered at $\infty$ determined by the distance.

The enclosure method was introduced in \cite{I1, IE} firstly for
inverse boundary value problems governed by elliptic equations.
In \cite{I4} the author has found its applications in several
inverse initial boundary value problems governed by heat and wave
equations in one-space dimensional case. Quite recently in
\cite{IK1, IK2} the method has been extended to three-dimensional
problems for the heat equation and in \cite{IE2} the framework for
inverse initial boundary value problem governed by the heat
equation together with a result when the background medium is
isotropic and inhomogeneous has been established.

This paper is a continuation of \cite{IE0} and the purpose
consists of two parts. First we extend the range of applications
of this new method to more general or another kind of obstacles.
We consider obstacles with a {\it dissipative boundary condition}
which is a generalization of the sound-hard obstacles; obstacles
with a {\it finite refractive index}, so-called, transparent
obstacles. Second we consider so-called {\it inverse
back-scattering problems} for those obstacles.  Therein the
observation data are measured in a limited time at the same
position as the initial data.  Needless to say this is a very
important class of inverse problems for waves and the solution has
many possibilities of applications, for example, nondestructive
testing such as diffraction tomography, subsurface radar,
microwave tomography, ocean acoustics tomography, etc. We show
rigorously that the method works also for inverse back-scattering
problems for obstacles mentioned above.

\subsection{Obstacle with a dissipative boundary condition}

First we consider an inverse obstacle scattering problem which is described by the classical wave equation outside an obstacle
with a dissipative boundary condition.

Let $D$ be a nonempty bounded open subset of $\Bbb R^3$ with smooth boundary such that
$\Bbb R^3\setminus\overline D$
is connected. Let $\gamma$ and $\beta$ be functions belonging to $L^{\infty}(\partial D)$ and satisfy $\gamma\ge 0$.
Let $0<T<\infty$.
Given $f\in L^2(\Bbb R^3)$ with compact support satisfying $\text{supp}\,f\cap\overline D=\emptyset$
let $u=u_f(x,t)$ denote the weak solution of the following initial boundary value problem for the classical wave equation:
$$\begin{array}{c}
\displaystyle
\partial_t^2u-\triangle u=0\,\,\text{in}\,(\Bbb R^3\setminus\overline D)\times\,]0,\,T[,\\
\\
\displaystyle
u(x,0)=0\,\,\text{in}\,\Bbb R^3\setminus\overline D,\\
\\
\displaystyle
\partial_tu(x,0)=f(x)\,\,\text{in}\,\Bbb R^3\setminus\overline D,\\
\\
\displaystyle
\frac{\partial u}{\partial\nu}-\gamma(x)\partial_tu-\beta(x)u=0\,\,\text{on}\,\partial D\times\,]0,\,T[.
\end{array}
\tag {1.1}
$$
Here $\nu$ denotes the outward normal to $D$ on $\partial D$.  Note that the boundary condition with $\beta\ge 0$ is called
the dissipative boundary condition, however, in this paper we include
also the case $\beta<0$ in this terminology.
Note also that a formal computation yields
$$\displaystyle
E'(t)=-2\int_{\partial D}\gamma(x)\vert\partial_t u\vert^2dS\le 0,
$$
where
$$\displaystyle
E(t)=\int_{\Bbb R^3\setminus\overline D}(\vert\partial_t u\vert^2+\vert\nabla u\vert^2)dx+\int_{\partial D}\beta(x)\vert u\vert^2dS,\,t\in [0,\,T].
$$

Let $\Omega$ be a bounded domain of $\Bbb R^3$ with a smooth boundary such that $\overline D\subset\Omega$ and $\Bbb R^3\setminus\overline\Omega$ is connected.
In this subsection we consider the following problem:

{\bf\noindent Inverse Problem I.}  Assume that $D$, $\gamma$ and $\beta$ on $\partial D$
are {\it unknown}.  Extract information about the location and shape of $D$ from the wave field
$u_f(x,t)$ given at all $x\in\partial\Omega$ and $t\in\,]0,\,T[$ for a fixed $f\in L^2(\Bbb R^3)$ with compact support
satisfying $\text{supp}\,f\cap\overline\Omega=\emptyset$.

Thus $\partial\Omega$ is a mathematical model of the closed
surface on which the wave field is measured.

Let $B$ denote an open ball.
In this paper as done in \cite{IE} we introduce two conditions on initial data $f$:

(I1)  $f(x)=0$ a.e. $x\in\Bbb R^3\setminus B$;

(I2)  there exists a positive constant $C$ such that $f(x)\ge C$ a. e. $x\in B$ or
$-f(x)\ge C$ a.e. $x\in B$.

Note that (I1) and (I2) imply that $\text{supp}\,f=\overline B$.

Set
$$\displaystyle
w(x)=w_f(x,\tau)=\int_0^Te^{-\tau t}u_f(x,t)dt,\,x\in\Bbb R^3\setminus\overline D,\,\tau>0.
\tag {1.2}
$$
We denote by $\nu$ also the outward normal to $\Omega$ on $\partial\Omega$.

We assume that

$\bullet$  there exists a positive constant $C'$
such that the one of the following two conditions
is satisfied:

(A1)  $\gamma(x)\le 1-C'$ a.e. $x\in\partial D$;

(A2)  $\gamma(x)\ge 1+C'$ a.e. $x\in\partial D$.

The first result of this subsection
is the following formula which employs $w$ together with $\partial w/\partial\nu$
on $\partial\Omega$.

\proclaim{\noindent Theorem 1.1.}
Let $B$ satisfy $\overline B\cap\overline\Omega=\emptyset$ and $f\in L^2(\Bbb R^3)$
satisfy both (I1) and (I2).
Let $v=v_f(\,\cdot\,,\tau)\in H^1(\Bbb R^3)$ denote the weak solution of the modified Helmholtz equation
$$\displaystyle
(\triangle-\tau^2)v+f=0\,\,\text{in}\,\Bbb R^3.
\tag {1.3}
$$

If the observation time $T$ satisfies
$$\displaystyle
T>2\text{dist}\,(D,B)-\text{dist}\,(\Omega,B),
\tag {1.4}
$$
then we have:

if (A1) is satisfied, then there exists a $\tau_0>0$ such that, for all $\tau\ge\tau_0$
$$\displaystyle
\int_{\partial\Omega}\left(\frac{\partial v}{\partial\nu}w
-\frac{\partial w}{\partial\nu}v\right)dS>0;
\tag {1.5}
$$

if (A2) is satisfied, then
then there exists a $\tau_0>0$ such that, for all $\tau\ge\tau_0$
$$\displaystyle
\int_{\partial\Omega}\left(\frac{\partial v}{\partial\nu}w
-\frac{\partial w}{\partial\nu}v\right)dS<0;
\tag {1.6}
$$

In both cases the formula
$$
\lim_{\tau\longrightarrow\infty}\frac{1}{2\tau}\log
\left\vert\int_{\partial\Omega}\left(\frac{\partial v}{\partial\nu}w
-\frac{\partial w}{\partial\nu}v\right)dS\right\vert
=-\text{dist}\,(D,B),
\tag {1.7}
$$
is valid.

\endproclaim

{\bf\noindent Remark 1.1.} $v$ in Theorem 1.1 is explicitly given by the formula
$$\displaystyle
v(x)=v_f(x,\tau)=\frac{1}{4\pi}\int_B\frac{e^{-\tau\vert x-y\vert}}{\vert x-y\vert}f(y)dy
\tag {1.8}
$$
and thus computable from $f$ in principle.

From the data $u_f(x,t)$ on $\partial\Omega\times\,]0,\,T[$ one
can compute $w$ together with $\partial w/\partial\nu$ on
$\partial\Omega$.
The computation process makes use of the uniqueness
of an initial boundary value problem for the wave equation in $(\Bbb R^3\setminus\overline\Omega)\times]0,\,T[$.
We do not repeat this description in detail and
see \cite{IE0} for the procedure.
Thus formula (1.7) is: a solution
to Inverse Problem I provided $\gamma$ satisfies (A1) or (A2);
an extension of Theorem 1.1 in \cite{IE0} which treated {\it sound-hard obstacles},
that is, the case when $\gamma=\beta=0$ on $\partial D$.

Note that $\text{dist}\,(D,\,B)+\sqrt{\vert\partial B\vert/4\pi}$
coincides with the distance from the center of $B$ to $D$.  Thus
by choosing a suitable initial data which is independent of $D$,
from (1.7) one gets the distance of $D$ to a given point in $\Bbb
R^3\setminus\overline\Omega$.

By Proposition 1.1. in \cite{IE0} (1.4) ensures that $T>l(\partial B, \partial D, \partial\Omega)$, where
$$\displaystyle
l(\partial B,\partial D,\partial\Omega)
=\inf\{\vert x-y\vert+\vert y-z\vert\,\vert\,x\in\partial B,\,y\in\partial D,\,z\in\partial\Omega\}.
$$
Thus (1.4) does not contradicts the {\it finite propagation
property} of the signal governed by the wave equation. 
It should be emphasized that we never make use of the finite
propagation property of the signal governed by the wave equation
in any form and this is an interesting point of our method since
we could find restriction (1.4).

It should be pointed out also that Theorem 1.1 gives a characterization whether
$\gamma<1$ or $\gamma>1$ in terms of a solution of the wave equation over a finite time interval
provided $\gamma$ is {\it constant} and $\gamma\not=1$.

Our method covers also another very important case which employs
the data observed on the support of the initial data not on
$\partial\Omega$, that is, the back-scattering one. More precisely
we consider the following problem.

{\bf\noindent Inverse Problem I'.}  Assume that $D$, $\gamma$ and $\beta$ on $\partial D$
are unknown.  Let $B$ satisfy $\overline
B\cap\overline D=\emptyset$. Extract information about the
location and shape of $D$ from the wave field $u_f(x,t)$ given at
all $x\in B$ and $t\in\,]0,\,T[$ for a fixed $f\in L^2(\Bbb R^3)$
with support $\overline B$.

A new finding in the enclosure method what we want to emphasize is the following result.

\proclaim{\noindent Theorem 1.2.}
Let $B$ satisfy $\overline B\cap\overline D=\emptyset$ and $f\in L^2(\Bbb R^3)$
satisfy both (I1) and (I2).
Let $v=v_f(\,\cdot\,,\tau)\in H^1(\Bbb R^3)$ be the weak solution of (1.3).
If $T$ satisfies
$$\displaystyle
T>2\text{dist}\,(D,B),
\tag {1.9}
$$
then we have:

if (A1) is satisfied, then there exists a $\tau_0>0$ such that, for all $\tau\ge\tau_0$
$$\displaystyle
\int_{B}f(w-v)dx>0;
\tag {1.10}
$$

if (A2) is satisfied, then there exists a $\tau_0>0$ such that, for all $\tau\ge\tau_0$
$$\displaystyle
\int_{B}f(w-v)dx<0;
\tag {1.11}
$$

In both cases the formula
$$
\lim_{\tau\longrightarrow\infty}\frac{1}{2\tau}\log\left\vert\int_{B}f(w-v)dx\right\vert
=-\text{dist}\,(D,B),
\tag {1.12}
$$
is valid.

\endproclaim

In formula (1.12) we make use of $w(x,\tau)$ for $x\in B$ which can be computed directly from
$u_f(x,t)$ for $(x,t)\in B\times\,]0,\,T[$ via (1.2).
Since $\text{supp}\,f=\overline B$,
this is a {\it back-scattering data} over finite time interval $]0,\,T[$.
Thus (1.12) gives a solution to Inverse Problem I'.  Note also that $B$ can be arbitrary small.

Theorem 1.1 ensures that even the case when $T$ is not greater
than $2\text{dist}\,(D,\,B)$, however satisfies (1.4) one can
extract $\text{dist}\,(D,\,B)$ via formula (1.7).  This is an
expression of an advantage of making use of data $u_f$ on the
whole of $\partial\Omega$.  In back-scattering case, restriction
(1.9) on $T$ is quite natural expected or {\it optimal} one
because of the finite propagation property of the signal in the
wave phenomena, however, note that in the proof we never make use
of this property in any form.

In the framework of the Lax-Phiilips scattering theory \cite{LPD},
Majda \cite{Mo} considered the case when $\beta=0$ and $D$ is {\it strictly
convex}.  He clarified the leading term of the
scattering amplitude $s(\theta,\omega,\lambda)$ at {\it high frequency} which is the Fourier
transform of the scattering kernel $S(s,\theta,\omega)$.  
However, our data are different from his ones and in this paper
we never assume that $D$ is strictly convex.

\subsection{Obstacle with a finite refractive index}

Second we consider an inverse scattering problem for an obstacle with
a finite refractive index.

Let $0<T<\infty$.  Given $f\in L^2(\Bbb R^3)$ with compact support
let $u=u_f(x,t)$ denote the weak solution of the following initial boundary value problem:
$$\begin{array}{c}
\displaystyle
\alpha(x)\partial_t^2u-\triangle u=0\,\,\text{in}\,\Bbb R^3\times\,]0,\,T[,\\
\\
\displaystyle
u(x,0)=0\,\,\text{in}\,\Bbb R^3,\\
\\
\displaystyle
\partial_tu(x,0)=f(x)\,\,\text{in}\,\Bbb R^3,
\end{array}
\tag {1.13}
$$
where $\alpha$ is a function belonging to $L^{\infty}(\Bbb R^3)$ and
satisfies $\alpha(x)\ge C$ a.e. $x\in\Bbb R^3$ for a positive constant $C$.

We assume that

$\bullet$  there exists a nonempty bounded open set $D$ with a smooth boundary
such that $\alpha(x)=1$ a.e. $x\in\Bbb R^3\setminus D$;

$\bullet$  there exists a positive constant $C'$ such that the one of the following two conditions
is satisfied:

(B1)  $\alpha(x)\le 1-C'$ a.e. $x\in D$;

(B2)  $\alpha(x)\ge 1+C'$ a.e. $x\in D$.

Thus $\alpha$ has a jump across $\partial D$.  $D$ is a mathematical model of an obstacle
with a finite refractive index.

In this subsection we present a solution to the following problem.

{\bf\noindent Inverse Problem II.}
Let $\Omega$ be a bounded domain of $\Bbb R^3$ with a smooth boundary such that $\overline D\subset\Omega$ and $\Bbb R^3\setminus\overline\Omega$ is connected.
Assume that both $D$ and $\alpha$ in $D$ are {\it unknown}.
Extract information about the location and shape of $D$ from the wave field
$u_f(x,t)$ given at all $x\in\partial\Omega$ and $t\in\,]0,\,T[$ for a fixed $f\in L^2(\Bbb R^3)$ with compact support
satisfying $\text{supp}\,f\cap\overline\Omega=\emptyset$.

Set
$$\displaystyle
w(x)=w_f(x,\tau)=\int_0^Te^{-\tau t}u_f(x,t)dt,\,x\in\Bbb R^3,\,\tau>0.
\tag {1.14}
$$

\proclaim{\noindent Theorem 1.3.}
Let $B$ satisfy $\overline B\cap\overline\Omega=\emptyset$ and $f\in L^2(\Bbb R^3)$
satisfy both (I1) and (I2).
Let $v=v_f(\,\cdot\,,\tau)\in H^1(\Bbb R^3)$ be the weak solution of (1.3).

If $T$ satisfy (1.4), then, we have:

if (B1) is satisfied, then there exists a $\tau_0>0$ such that, for all $\tau\ge\tau_0$
$$\displaystyle
\int_{\partial\Omega}\left(\frac{\partial v}{\partial\nu}w
-\frac{\partial w}{\partial\nu}v\right)dS>0;
\tag {1.15}
$$

if (B2) is satisfied, then there exists a $\tau_0>0$ such that,
for all $\tau\ge\tau_0$
$$\displaystyle
\int_{\partial\Omega}\left(\frac{\partial v}{\partial\nu}w
-\frac{\partial w}{\partial\nu}v\right)dS<0.
\tag {1.16}
$$

In both cases we have
$$\displaystyle
\lim_{\tau\longrightarrow\infty}\frac{1}{2\tau}\log\left\vert\int_{\partial\Omega}\left(\frac{\partial v}{\partial\nu}w
-\frac{\partial w}{\partial\nu}v\right)dS\right\vert
=-\text{dist}\,(D,B).
\tag {1.17}
$$
\endproclaim

The procedure for the computation of both $w$ and $\partial w/\partial\nu$ on $\partial\Omega$ from
$u_f(x,t)$ given at all $x\in\partial\Omega$ and $t\in\,]0,\,T[$ is the same as Theorem 1.1.
Thus this is a solution to Inverse Problem II.

The method works also for the following problem.

{\bf\noindent Inverse Problem II'.}
Assume that both $D$ and $\alpha$ in $D$ are unknown.
Let $B$ satisfy $\overline B\cap\overline D=\emptyset$.
Extract information about the location and shape of $D$ from the wave field
$u_f(x,t)$ given at all $x\in B$ and $t\in\,]0,\,T[$ for a fixed $f\in L^2(\Bbb R^3)$ with support $\overline B$.

The following result gives a solution to Inverse Problem II'.

\proclaim{\noindent Theorem 1.4.}
Let $B$ satisfy $\overline B\cap\overline D=\emptyset$ and $f\in L^2(\Bbb R^3)$
satisfy both (I1) and (I2).
Let $v=v_f(\,\cdot\,,\tau)\in H^1(\Bbb R^3)$ be the weak solution of (1.3).

If $T$ satisfies (1.9), then we have:

if (B1) is satisfied, then there exists a $\tau_0>0$ such that, for all $\tau\ge\tau_0$
$$\displaystyle
\int_{B}f(w-v)dx>0;
\tag {1.18}
$$

if (B2) is satisfied, then there exists a $\tau_0>0$ such that,
for all $\tau\ge\tau_0$
$$\displaystyle
\int_{B}f(w-v)dx<0.
\tag {1.19}
$$

In both cases we have
$$\displaystyle
\lim_{\tau\longrightarrow\infty}\frac{1}{2\tau}\log
\left\vert\int_{B}f(w-v)dx\right\vert
=-\text{dist}\,(D,B).
\tag {1.20}
$$

\endproclaim

In the framework of the Lax-Phillips scattering theory there is a result by Majda-Taylor \cite{MT}
for the case when $\alpha$ is smooth on $\overline D$ and $\alpha(x)\not=1$.  However,
the data are given by the back-scattering kernel $S(s,-\omega,\omega)$, $\omega\in S^2$ and the obtained information is the value of the support function of unknown obstacles
and thus different from ours.
For the study of the leading term of the scattering amplitude at high frequency for a strictly convex obstacle with a finite refractive index
see Majda-Taylor \cite{MT} ($\alpha(x)>1$ on $\overline D$) and
Petkov \cite{P1,P2} ($\alpha(x)<1$ on $\overline D$).

\subsection{Construction of the paper}

A brief outline of this paper is as follows. Theorems 1.1 and 1.2
are proved in Section 2. In Subsection 2 first we formulate what
we mean by the weak solution of (1.1). It is based on the notion
of the weak solution in \cite{DL} and an application of the
existence theory therein.  We see that $w$ given by (1.2)
satisfies the modified Helmholtz equation with an {\it unknown}
inhomogeneous term outside $D$ and boundary data on $\partial D$
in a weak sense. In Subsection 2.2 two expressions for
the integral
$$\displaystyle
\int_{\partial\Omega}
\left(\frac{\partial v}{\partial\nu}w-\frac{\partial w}{\partial\nu}v\right)dS
\tag {1.21}
$$
are established.  One is called the {\it local expression} which yields a bound
of the absolute value of (1.21).
Another is called the {\it global expression} and yields the lower bound of (1.21)
in the case when (A1) is satisfied: upper bound of that
in the case when (A2) is satisfied.  Those bounds are presented
in Subsections 2.4 and 2.5, respectively.  Theorem 1.1 is a direct consequence of those bounds.
The proof of Theorem 1.2 is described in Subsection 3.6.  It is based on the following
asymptotic formula which connects two data in Theorems 1.1 and 1.2
for an arbitrary fixed $T<\infty$:
$$\displaystyle
\int_{\Bbb R^3\setminus\overline D}f(w-v)dx=
\int_{\partial\Omega}
\left(\frac{\partial v}{\partial\nu}w-\frac{\partial w}{\partial\nu}v\right)dS
+O(\tau^{-1}e^{-\tau T}).
\tag {1.22}
$$
Using this together with the obtained bounds in Subsections 2.4 and 2.5,
one gets immediately the conclusion of Theorem 1.2.

Section 3 is devoted to the proof of Theorems 1.3 and 1.4. The
order of the presentation is parallel to that of Section 2.
Starting with the formulation of the weak solution of (1.13), we
present two expressions of (1.21) for $w$ given by (1.14).  Both
expressions are of global type since they involve integrals over
the whole space. Using those expressions, we give a bound of the
absolute value of(1.21) in Subsections 4.2; the upper/lower bound
of (1.21) for case (B1)/(B2).  Theorem 1.3 is a direct consequence
of those bounds and an asymptotic formula which connects two data
in Theorems 1.3 and 1.4 similar to (1.22) enables us to obtain all
the conclusions of Theorem 1.4.

Section 4 consists of some remarks and open problems.
In Subsection 4.1 we show that the case $\gamma=1$ is {\it exceptional}
in one-space dimensional case.
Moreover, the complete asymptotic expansion of (1.21) in one
space-dimensional case is presented. In Subsection 4.2
corresponding to results in one-space dimensional case we propose
some open problems and point out related results in some
references together with future direction of our method.

\section{Proof of Theorems 1.1 and 1.2}

\subsection{The weak solution and the governing equation for $w_f$}

The contents of this subsection almost parallel to the corresponding
parts of Section 2 in \cite{IE0} which is the case when $\gamma=\beta=0$
on $\partial D$.

We write $u'$ instead of $\partial_tu$.
It has been observed in \cite{IE0} that
for the enclosure method the notion of the weak solution
described on pp 552-66 in \cite{DL} is appropriate since we consider only the solution
on a {\it finite time interval}.

Set $V=H^1(\Bbb R^3\setminus\overline D)$ and $H=L^2(\Bbb
R^3\setminus\overline D)$. Applying Theorem 1 on p 558 in
\cite{DL} to (1.1), we know that: given $u^0\in V$ and $u^1\in H$
there exists a unique $u$ satisfying
$$
\displaystyle
u\in L^2(0,\,\,T;V),\,\,
u'\in L^2(0,\,\,T;V),\,\,
u''\in L^2(0,\,\,T;V' )
$$
such that, for all $\phi\in V$
$$\displaystyle
<u''(t),\phi>
+a(u(t),\phi)+b(u'(t),\phi)=0,\,\,\text{a.e.}\,t\in]0,\,\,T[,
\tag {2.1}
$$
and $u(0)=u^0$ and $u'(0)=u^1$, where
$$\displaystyle
a(u,v)=\int_{\Bbb R^3\setminus\overline D}\nabla u\cdot\nabla vdx
+\int_{\partial D}\beta uvdS,\,\,u,v\in V
$$
and
$$\displaystyle
b(u,v)=\int_{\partial D}\gamma uvdS,\,\,u,v\in V.
$$
Note that this is the case when $f=0$ in their notation
and from their proof one knows that
the restriction (5.8) on p 553 for $b$ is {\it redundant} in this case, that means that
$b_0$ therein can be identically zero and $b=b_1$.  Note that at this step we do not make use
of condition $\gamma\ge 0$ on $\partial D$.

In this section we say that this $u=u_f$ for $u^0=0$ and
$u^1=f$ is the weak solution of (1.1).
Then $w=w_f$ given by (1.2) belongs to $V$; it follows from integration by parts and (2.1)
multiplied by $e^{-\tau t}$ that $w$ satisfies for all $\phi\in V$
$$\begin{array}{c}
\displaystyle
\int_{\Bbb R^3\setminus\overline D}
\nabla w\cdot\nabla\phi dx
+\int_{\partial D}c(x,\tau)w\phi dS
+\int_{\Bbb R^3\setminus\overline D}(\tau^2w-f)\phi dx\\
\\
\displaystyle
=-e^{-\tau T}\int_{\Bbb R^3\setminus\overline D}F\phi dx
-e^{-\tau T}\int_{\partial D}G\phi dS,
\end{array}
\tag {2.2}
$$
where
$$\begin{array}{c}
\displaystyle
c(x,\tau)=\gamma(x)\tau+\beta(x),\\
\\
\displaystyle
F(x,\tau)=u'(x,T)+\tau u(x,T),\\
\\
\displaystyle
G(x)=\gamma(x)u(x,T).
\end{array}
$$
This means that, in a weak sense,  $w$ satisfies
$$\begin{array}{c}
\displaystyle
(\triangle-\tau^2)w+f=e^{-\tau T}F(x,\tau)\,\,
\text{in}\,\Bbb R^3\setminus\overline D,\\
\\
\displaystyle
\frac{\partial w}{\partial\nu}
=c(x,\tau)w+e^{-\tau T}G(x)\,\,\text{on}\,\partial D.
\end{array}
$$
In particular, from (2.2) for $\phi\in C_0^{\infty}(\Bbb R^3\setminus\overline D)$,
we have $(\triangle-\tau^2)w+f(x)=e^{-\tau T}F(x,\tau)$ in $\Bbb R^3\setminus\overline D$
in the sense of distribution and hence $\triangle w\in L^2(\Bbb R^3\setminus\overline D)$.
This yields $w\in H^2_{\text{loc}}(\Bbb R^3\setminus\overline D)$ and
$(\triangle-\tau^2)w+f(x)=e^{-\tau T}F(x,\tau)$ a.e. $x\in\Bbb R^3\setminus\overline D$.
Thus we can define $\partial w/\partial\nu\vert_{\partial\Omega}$ as
$\nabla w\vert_{\partial\Omega}\cdot\nu\in H^{1/2}(\partial\Omega)$,
where $\nabla w\vert_{\partial\Omega}$ is the trace of $\nabla w$ onto $\partial\Omega$
(\cite{G}).

\subsection{Local and global expressions of (1.21)}

In the following proposition we do not assume that $\text{supp}\,f\cap\overline D=\emptyset$ nor
$\text{supp}\,f\cap\overline\Omega=\emptyset$.

\proclaim{\noindent Proposition 2.1.}
We have two expressions:
$$
\begin{array}{c}
\displaystyle
\int_{\partial\Omega}\left(\frac{\partial v}{\partial\nu}w
-\frac{\partial w}{\partial\nu}v\right)dS
\\
\\
\displaystyle
=\int_{\partial D}\left(\frac{\partial v}{\partial\nu}-c(x,\tau)v\right)wdS
-\int_{\Omega\setminus\overline D}(w-v)fdx\\
\\
\displaystyle
-e^{-\tau T}\int_{\partial D}GvdS-e^{-\tau T}\int_{\Omega\setminus\overline D}Fvdx;
\end{array}
\tag {2.3}
$$
$$\begin{array}{c}
\displaystyle
\int_{\partial\Omega}\left(\frac{\partial v}{\partial\nu}w
-\frac{\partial w}{\partial\nu}v\right)dS\\
\\
\displaystyle
=\int_D(\vert\nabla v\vert^2+\tau^2\vert v\vert^2)dx-\int_{\partial D}c(x,\tau)\vert v\vert^2dS
-e^{-\tau T}\int_{\partial D}GvdS
\\
\\
\displaystyle
+\int_{\Bbb R^3\setminus\overline D}\left\{\vert\nabla(w-v)\vert^2
+\tau^2\vert w-v\vert^2\right\}dx
+\int_{\partial D}c(x,\tau)\vert w-v\vert^2dS\\
\\
\displaystyle
+e^{-\tau T}\int_{\Bbb R^3\setminus\overline D}F(w-v)dx-e^{-\tau T}\int_{\Omega\setminus\overline D}Fvdx+
e^{-\tau T}\int_{\partial D}G(w-v)dS
\\
\\
\displaystyle
-\int_{\Omega\setminus\overline D}(w-v)fdx-\int_Dfvdx.
\end{array}
\tag {2.4}
$$

\endproclaim

{\it\noindent Proof.}
First we give a proof of (2.3).
Let $\varphi\in H^1(\Omega\setminus\overline D)$ satisfy $\varphi=0$ on $\partial\Omega$
in the sense of the trace.  Since the zero extension to $\Bbb R^3\setminus\overline D$
of this $\varphi$ belongs to $V$ (\cite{G}), we have from (2.2)
$$\begin{array}{c}
\displaystyle
0=\int_{\Omega\setminus\overline D}
\nabla w\cdot\nabla\varphi dx
+\int_{\partial D}c(x,\tau)w\varphi dS\\
\\
\displaystyle
+\int_{\Omega\setminus\overline D}(\tau^2w-f+e^{-\tau T}F)\varphi dx
+e^{-\tau T}\int_{\partial D}G\varphi dS.
\end{array}
\tag {2.5}
$$
Choose $\chi\in C^{\infty}_0(\Bbb R^3)$ such that $\chi=1$ in a neighbourhood
of $\partial\Omega$ and $\chi(x)\equiv 0$ in a neighbourhood of $\overline D$.
Since $(1-\chi)v\vert_{\Omega\setminus\overline D}$
vanishes in a neighbourhood of $\partial\Omega$, (2.5) for
$\varphi=(1-\chi)v\vert_{\Omega\setminus\overline D}$ is valid.

On the other hand, since $\chi v$ vanishes in a neighbourhood of $\overline D$
and $w\in H^2_{\text{loc}}(\Bbb R^3\setminus\overline D)$,
integration by parts yields
$$\begin{array}{c}
\displaystyle
\int_{\Omega\setminus\overline D}\nabla w\cdot\nabla(\chi v)dx
=\int_{\partial\Omega}\frac{\partial w}{\partial\nu}\chi vdS
-\int_{\Omega\setminus\overline D}(\triangle w)\chi vdx\\
\\
\displaystyle
=\int_{\partial\Omega}\frac{\partial w}{\partial\nu}vdS
-\int_{\Omega\setminus\overline D}
(\tau^2w-f+e^{-\tau T}F)\chi vdx,
\end{array}
$$
that is
$$\displaystyle
\int_{\partial\Omega}\frac{\partial w}{\partial\nu}vdS
=\int_{\Omega\setminus\overline D}\nabla w\cdot\nabla(\chi v)dx
+\int_{\Omega\setminus\overline D}
(\tau^2w-f+e^{-\tau T}F)\chi vdx.
$$
From this and (2.5) for $\varphi=(1-\chi)v\vert_{\Omega\setminus\overline D}$
which satisfies $\varphi=v$ on $\partial D$ we obtain
$$\begin{array}{c}
\displaystyle
\int_{\partial\Omega}
\frac{\partial w}{\partial\nu}vdS
=\int_{\Omega\setminus\overline D}\nabla w\cdot\nabla vdx
+\int_{\partial D}c(x,\tau)wvdS\\
\\
\displaystyle
+\int_{\Omega\setminus\overline D}(\tau^2w-f+e^{-\tau T}F)vdx
+\int_{\partial D}e^{-\tau T}GvdS.
\end{array}
\tag {2.6}
$$

By the trace theorem one can choose $\tilde{w}\in H^1(\Omega)$ such that $\tilde{w}=w$ in $\Omega\setminus\overline D$.
Note that $\tilde{w}=w$ on $\partial\Omega$ and $\partial D$ in the sense of the trace.
Since $v\in H^2(\Omega)$, we have
$$\begin{array}{c}
\displaystyle
\int_{\partial\Omega}\frac{\partial v}{\partial\nu} wdS
=\int_{\partial\Omega}\frac{\partial v}{\partial\nu} \tilde{w}dS\\
\\
\displaystyle
=\int_{\Omega}\triangle v \tilde{w}dx+\int_{\Omega}\nabla v\cdot\nabla\tilde{w}dx\\
\\
\displaystyle
=\int_{\Omega\setminus\overline D}(\tau^2 v-f)wdx+\int_{\Omega\setminus\overline D}\nabla v\cdot\nabla wdx\\
\\
\displaystyle
+\int_{D}(\tau^2 v-f)\tilde{w}dx+\int_{D}\nabla v\cdot\nabla\tilde{w}dx.
\end{array}
$$
On the other hand we have
$$\begin{array}{c}
\displaystyle
\int_D\nabla v\cdot\nabla\tilde{w}dx
=\int_{\partial D}
\frac{\partial v}{\partial\nu}\tilde{w}dS-\int_D(\triangle v)\tilde{w}dx
\\
\\
\displaystyle
=\int_{\partial D}
\frac{\partial v}{\partial\nu}wdS-\int_D(\tau^2 v-f)\tilde{w}dx,
\end{array}
$$
that is,
$$\displaystyle
\int_{\partial D}\frac{\partial v}{\partial\nu}wdS
=\int_D(\tau^2 v-f)\tilde{w}dx
+\int_D\nabla v\cdot\nabla\tilde{w}dx.
$$
Therefore, we obtain
$$\displaystyle
\int_{\partial\Omega}\frac{\partial v}{\partial\nu} wdS
=\int_{\Omega\setminus\overline D}(\tau^2 v-f)wdx+\int_{\Omega\setminus\overline D}\nabla v\cdot\nabla wdx
+\int_{\partial D}\frac{\partial v}{\partial\nu}wdS.
\tag {2.7}
$$
A combination of (2.6) and (2.7) gives (2.3).

Next we give a proof of (2.4).
Write
$$\displaystyle
\int_{\partial D}\frac{\partial v}{\partial\nu}wdS
=\int_{\partial D}\frac{\partial v}{\partial\nu}(w-v)dS+\int_{\partial D}\frac{\partial v}{\partial\nu}vdS.
\tag {2.8}
$$
By the trace theorem given $\phi\in V$ one can choose $\tilde{\phi}\in H^1(\Bbb R^3)$ such that $\phi=\tilde{\phi}\vert_{\Bbb R^3\setminus\overline D}$.
Since $v\in H^2(D)$ and $\triangle v-\tau^2v+f=0$ a.e. $x\in D$ and $\phi=\tilde{\phi}$ on $\partial D$ in the sense of the trace,
integration by parts yields
$$\displaystyle
\int_{\partial D}\frac{\partial v}{\partial\nu}\phi dS
=\int_D(\nabla v\cdot\nabla\tilde{\phi}dx+\tau^2v\tilde{\phi})dx-\int_Df\tilde{\phi}dx.
\tag {2.9}
$$
On the other hand $v$ satisfies
$$\displaystyle
-\int_{\Bbb R^3}\nabla v\cdot\nabla\varphi dx-\tau^2\int_{\Bbb R^3}v\varphi dx=-\int_{\Bbb R^3}f\varphi dx,\,\,\forall\varphi\in H^1(\Bbb R^3).
$$
Substituting $\varphi=\tilde{\phi}$ into this identity and dividing $\Bbb R^3=D\cup(\Bbb R^3\setminus D)$, we obtain
$$\displaystyle
\int_{\Bbb R^3\setminus\overline D}(\nabla v\cdot\nabla\phi+\tau^2 v\phi)dx=\int_{\Bbb R^3\setminus\overline D}
f\phi dx-\int_{D}(\nabla v\cdot\nabla\tilde{\phi}+\tau^2 v\tilde{\phi})dx
+\int_{D}f\tilde{\phi} dx.
$$
A combination of this and (2.9) gives, for all $\phi\in V$,
$$\displaystyle
\int_{\partial D}\frac{\partial v}{\partial\nu}\phi dS
+\int_{\Bbb R^3\setminus\overline D}(\nabla v\cdot\nabla\phi
+\tau^2v\phi)dx=\int_{\Bbb R^3\setminus\overline D}f\phi dx.
\tag {2.10}
$$
Combining this with (2.2), we obtain
$$\begin{array}{c}
\displaystyle
\int_{\partial D}\left(\frac{\partial v}{\partial\nu}-c(x,\tau)v\right)\phi dS
=\int_{\partial D}c(x,\tau)\epsilon\phi dS
+\int_{\Bbb R^3\setminus\overline D}\left\{\nabla\epsilon\cdot\nabla\phi+\tau^2\epsilon
\phi\right\}dx\\
\\
\displaystyle
+e^{-\tau T}\int_{\Bbb R^3\setminus\overline D}F\phi dx
+e^{-\tau T}\int_{\partial D}G\phi dS,
\end{array}
\tag {2.11}
$$
where $\epsilon=w-v$.  This means that $\epsilon$ satisfies, in a weak sense
$$\begin{array}{c}
\displaystyle
(\triangle-\tau^2)\epsilon=e^{-\tau T}F\,\,\text{in}\,\Bbb R^3\setminus\overline D,
\\
\\
\displaystyle
\frac{\partial\epsilon}{\partial\nu}
-c(x,\tau)\epsilon=-\frac{\partial v}{\partial\nu}+c(x,\tau)v+e^{-\tau T}G
\,\,\text{on}\,\partial D.
\end{array}
$$
Substituting $\epsilon$ for $\phi$ in (2.11), we obtain
$$\begin{array}{c}
\displaystyle
\int_{\partial D}\frac{\partial v}{\partial\nu}\epsilon dS
=\int_{\partial D}c(x,\tau)(\vert\epsilon\vert^2+v\epsilon)dS+\int_{\Bbb R^3\setminus\overline D}(\vert\nabla\epsilon\vert^2
+\tau^2\vert\epsilon\vert^2)dx\\
\\
\displaystyle
+e^{-\tau T}\int_{\Bbb R^3\setminus\overline D}F\epsilon dx
+e^{-\tau T}\int_{\partial D}G\epsilon dS
\end{array}
\tag {2.12}
$$
Now from this together with (2,3), (2.8) and the identities
$v\epsilon=-\vert v\vert^2+vw$ and
$$\begin{array}{c}
\displaystyle
\int_{\partial D}\frac{\partial v}{\partial\nu}vdS
=\int_D(\vert\nabla v\vert^2+\tau^2\vert v\vert^2)dx-\int_Dfvdx,
\end{array}
\tag {2.13}
$$
we obtain (2.4).

\noindent
$\Box$

Now assume that $\text{supp}\,f\cap\overline\Omega=\emptyset$. From (2.3) and (2.4) we have the following two expressions which we call the {\it local expression}
and {\it global expression}, respectively:
$$\begin{array}{c}
\displaystyle
\int_{\partial\Omega}\left(\frac{\partial v}{\partial\nu}w-\frac{\partial w}{\partial\nu}v\right)dS\\
\\
\displaystyle
=\int_{\partial D}\left(\frac{\partial v}{\partial\nu}-c(x,\tau)v\right)wdS
-e^{-\tau T}\int_{\partial D}GvdS-e^{-\tau T}\int_{\Omega\setminus\overline D}Fvdx;
\end{array}
\tag {2.14}
$$
$$\begin{array}{c}
\displaystyle
\int_{\partial\Omega}\left(\frac{\partial v}{\partial\nu}w
-\frac{\partial w}{\partial\nu}v\right)dS\\
\\
\displaystyle
=\int_D(\vert\nabla v\vert^2+\tau^2\vert v\vert^2)dx-\int_{\partial D}c(x,\tau)\vert v\vert^2dS-e^{-\tau T}\int_{\partial D}GvdS\\
\\
\displaystyle
+\int_{\Bbb R^3\setminus\overline D}\left\{\vert\nabla(w-v)\vert^2
+\tau^2\vert w-v\vert^2\right\}dx
+\int_{\partial D}c(x,\tau)\vert w-v\vert^2dS\\
\\
\displaystyle
+e^{-\tau T}\int_{\Bbb R^3\setminus\overline D}F(w-v)dx-e^{-\tau T}\int_{\Omega\setminus\overline D}Fvdx
+e^{-\tau T}\int_{\partial D}G(w-v)dS.
\end{array}
\tag {2.15}
$$

For convenience we set
$$\displaystyle
I(\tau)=\int_{\partial\Omega}\left(\frac{\partial v}{\partial\nu}w
-\frac{\partial w}{\partial\nu}v\right)dS.
$$
We make use of (2.14) to give an estimation of $\vert I(\tau)\vert$ from above
and (2.15) of $I(\tau)$ from below/above when (A1)/(A2) is satisfied.
This is the role of (2.14) and (2.15).

\subsection{An estimate of $\vert I(\tau)\vert$ from above}

In this subsection we derive the following estimate as $\tau\longrightarrow\infty$:
$$\displaystyle
e^{2\tau\text{dist}\,(D,\,B)}\vert I(\tau)\vert=O(\tau^{3/2}).
\tag {2.16}
$$

\proclaim{\noindent Lemma 2.1.}
Let $\epsilon=w-v$.
As $\tau\longrightarrow\infty$ we have
$$\displaystyle
\Vert\epsilon\Vert_{L^2(\Bbb R^3\setminus\overline D)}^2
=O(e^{-2\tau T}+e^{-2\tau\text{dist}\,(D,\,B)})
\tag {2.17}
$$
and
$$\displaystyle
\Vert\nabla\epsilon\Vert_{L^2(\Bbb R^3\setminus\overline D)}^2
=O(\tau^2(e^{-2\tau T}+e^{-2\tau\text{dist}\,(D,\,B)})).
\tag {2.18}
$$

\endproclaim

{\it\noindent Proof.}
Rewrite (1.12) as
$$\begin{array}{c}
\displaystyle
\int_{\Bbb R^3\setminus\overline D}(\vert\nabla\epsilon\vert^2
+\tau^2\vert\epsilon\vert^2)dx
=-\int_{\partial D}c(x,\tau)(\vert\epsilon\vert^2+v\epsilon)dS
+\int_{\partial D}\frac{\partial v}{\partial\nu}\epsilon dS
\\\\
\displaystyle
-e^{-\tau T}\left(\int_{\Bbb R^3\setminus\overline D}F\epsilon dx+\int_{\partial D}G\epsilon dS\right).
\end{array}
\tag {2.19}
$$
We have
$$\begin{array}{c}
\displaystyle
-\int_{\partial D}c(x,\tau)(\vert\epsilon\vert^2+v\epsilon)dS
+\int_{\partial D}\frac{\partial v}{\partial\nu}\epsilon dS\\
\\
\displaystyle
=-\tau\int_{\partial D}\gamma(\vert\epsilon\vert^2+v\epsilon)dS
+\int_{\partial D}\left(\frac{\partial v}{\partial\nu}-\beta\vert\epsilon\vert^2-\beta v\epsilon\right)dS.
\end{array}
\tag {2.20}
$$
Completing the square and $\gamma\ge 0$ on
$\partial D$ give
$$\begin{array}{c}
\displaystyle
-\int_{\partial D}\gamma(\vert\epsilon\vert^2+v\epsilon)dS
=-\int_{\partial D}\gamma\left\vert \epsilon+\frac{v}{2}\right\vert^2dS
+\int_{\partial D}\frac{\gamma}{4}\vert v\vert^2dS\\
\\
\displaystyle
\le\int_{\partial D}\frac{\gamma}{4}\vert v\vert^2dS.
\end{array}
\tag {2.21}
$$
Moreover we have
$$\displaystyle
\int_{\partial D}\left(\frac{\partial v}{\partial\nu}\epsilon-\beta\vert\epsilon\vert^2
-\beta v\epsilon\right)dS
\le C\int_{\partial D}\left(\vert\epsilon\vert^2+\vert v\vert^2
+\left\vert\frac{\partial v}{\partial\nu}\right\vert^2\right)dS,
\tag {2.22}
$$
$$\displaystyle
e^{-\tau T}\left\vert\int_{\partial D} G\epsilon dS\right\vert
\le\frac{1}{2}\left(\int_{\partial D}\vert\epsilon\vert^2 dS+e^{-2\tau T}\int_{\partial D}\vert G\vert^2dS\right)
\tag {2.23}
$$
and
$$\displaystyle
e^{-\tau T}\left\vert\int_{\Bbb R^3\setminus\overline D} F\epsilon dx\right\vert
\le\frac{1}{2}\left(\eta^2\int_{\Bbb R^3\setminus\overline D}\vert\epsilon\vert^2 dx+e^{-2\tau T}\eta^{-2}\int_{\Bbb R^3\setminus\overline D}\vert F\vert^2dx\right),\,\,\forall\eta>0.
\tag {2.24}
$$

From (2.19)-(2.24) we obtain
$$\begin{array}{c}
\displaystyle
\int_{\Bbb R^3\setminus\overline D}\vert\nabla\epsilon\vert^2dx+\left(\tau^2-\frac{\eta^2}{2}\right)\int_{\Bbb R^3\setminus\overline D}\vert\epsilon\vert^2dx\\
\\
\displaystyle
\le
\left(C+\frac{1}{2}\right)\int_{\partial D}\vert\epsilon\vert^2dS+
\int_{\partial D}\left(C+\frac{\tau\gamma}{4}\right)\vert v\vert^2dS
+C\int_{\partial D}\left\vert\frac{\partial v}{\partial\nu}\right\vert^2dS
\\
\\
\displaystyle
+\frac{e^{-2\tau T}}{2}\left(\int_{\partial D}\vert G\vert^2dS
+\frac{1}{2\eta^2}\int_{\Bbb R^3\setminus\overline D}\vert F\vert^2dx\right).
\end{array}
\tag {2.25}
$$
Here we cite the well known inequality (e.g., \cite{G}) that there
exists a positive constant $K=K(\Omega\setminus\overline D)$ such
that, for all $z\in H^1(\Omega\setminus\overline D)$
$$\displaystyle
\int_{\partial D}\vert z\vert^2dS
\le K\left(\eta^2\int_{\Omega\setminus\overline D}\vert\nabla z\vert^2dx
+\eta^{-2}\int_{\Omega\setminus\overline D}\vert z\vert^2 dx\right),\,\,\forall\eta>0.
\tag {2.26}
$$
Applying (2.26) to the first term in the right-hand side of (2.25), we have
$$
\begin{array}{c}
\displaystyle
\int_{\Bbb R^3\setminus\overline\Omega}\vert\nabla\epsilon\vert^2dx+
\left(\tau^2-\frac{\eta^2}{2}\right)\int_{\Bbb R^3\setminus\overline\Omega}\vert\epsilon\vert^2dx\\
\\
\displaystyle
+
\left\{1-K\eta^2\left(C+\frac{1}{2}\right)\right\}
\int_{\Omega\setminus\overline D}\vert\nabla\epsilon\vert^2dx+
\left\{\left(\tau^2-\frac{\eta^2}{2}\right)-K\eta^{-2}\left(C+\frac{1}{2}\right)\right\}
\int_{\Omega\setminus\overline D}\vert\epsilon\vert^2dx\\
\\
\displaystyle
\le
\int_{\partial D}\left\{\left(C+\frac{\tau\gamma}{4}\right)\vert v\vert^2
+C\left\vert\frac{\partial v}{\partial\nu}\right\vert^2\right\}dS
+\frac{e^{-2\tau T}}{2}\left(\int_{\partial D}\vert G\vert^2dS
+\frac{1}{2\eta^2}\int_{\Bbb R^3\setminus\overline D}\vert F\vert^2dx\right).
\end{array}
\tag {2.27}
$$
Now choosing a small $\eta$, from (1.8) and (2.27) we obtain (2.17) and (2.18) as $\tau\longrightarrow\infty$.

\noindent
$\Box$

Now from (2.26) for $\eta=1/\sqrt{\tau}$, (2.17) and (2.18) we obtain
$$\displaystyle
\Vert \epsilon\Vert_{L^2(\partial D)}^2
=O(\tau(e^{-2\tau T}+e^{-2\tau\text{dist}\,(D,\,B)})).
\tag {2.28}
$$
This together with $\Vert v\Vert_{L^2(\partial D)}=O(e^{-\tau\text{dist}\,(D,\,B)})$ yields
$$\displaystyle
\Vert w\Vert_{L^2(\partial D)}^2
=O(\tau(e^{-2\tau T}+e^{-2\tau\text{dist}\,(D,\,B)}))
\tag {2.29}
$$
Since
$$\displaystyle
\left\Vert\frac{\partial v}{\partial\nu}\right\Vert_{L^2(\partial D)}=O(\tau e^{-\tau\text{dist}\,(D,B)}),
$$
it follows from (2.29) that
$$\displaystyle
\left\vert
\int_{\partial D}\frac{\partial v}{\partial\nu}wdS\right\vert
+\left\vert
\int_{\partial D}c(x,\tau)vwdS\right\vert
=O(\tau^{3/2}(e^{-2\tau\text{dist}\,(D,\,B)}+e^{-\tau(T+\text{dist}\,(D,\,B)})).
\tag {2.30}
$$
Moreover we have
$$\displaystyle
e^{-\tau T} \left\vert\int_{\partial D} GvdS\right\vert
=O(e^{-\tau(T+\text{dist}\,(D,\,B))})
\tag {2.31}
$$
and since
$$\displaystyle
\left\Vert v\right\Vert_{L^2(\Omega\setminus\overline D)}
=O(e^{-\tau\text{dist}\,(\Omega,\,B)}),
$$
we obtain
$$\displaystyle
e^{-\tau T}
\left\vert\int_{\Omega\setminus\overline D}
Fvdx\right\vert
=O(\tau e^{-\tau(T+\text{dist}\,(\Omega,\,B))}).
\tag {2.32}
$$

Now from (2.14), (2.30), (2.31) and (2.32) we obtain
$$\displaystyle
\vert I(\tau)\vert
=O(\tau^{3/2}(e^{-2\tau\text{dist}\,(D,B)}+e^{-\tau(T+\text{dist}\,(D,\,B))})
+\tau e^{-\tau(T+\text{dist}\,(\Omega,\,B))})
$$
and thus
$$\displaystyle
e^{2\tau\text{dist}\,(D,\,B)}\vert I(\tau)\vert
=O(\tau^{3/2}(1+ e^{-\tau(T-\text{dist}\,(D,\,B))})
+\tau e^{-\tau(T-2\text{dist}\,(D,\,B)+\text{dist}\,(\Omega,\,B))}).
\tag {2.33}
$$
Write
$$\displaystyle
T-\text{dist}\,(D,\,B)
=(T-2\text{dist}\,(D,\,B)+\text{dist}\,(\Omega,\,B))
+(\text{dist}\,(D,\,B)-\text{dist}\,(\Omega,\,B)).
$$
This together with $\text{dist}\,(D,\,B)>\text{dist}\,(\Omega,\,B)$ and (1.4) gives
$T>\text{dist}\,(D,\,B)$.
Now from this together with (1.4) and (2.33) we obtain (2.16).

{\bf\noindent Remark 2.1.}
From this last part we know that, for the proof of (2.16) it suffice to assume only
$T\ge 2\text{dist}\,(D,\,B)-\text{dist}\,(\Omega,\,B)$ instead of stricter condition
(1.4).

\subsection{Case (A1).  An estimate of $I(\tau)$ from below}

In this subsection we prove that: there exists a positive constant $\mu$
such that
$$\displaystyle
\liminf_{\tau\longrightarrow\infty}\tau^{\mu}e^{2\tau\text{dist}\,(D,\,B)}I(\tau)>0
\tag {2.34}
$$
provided $\gamma$ satisfies (A1).

Having (2.16) and (2.34), we immediately obtain (1.5) and (1.7).

For the proof of (2.34) first we study the asymptotic behaviour of the following integral
$$\displaystyle
J(\tau)\equiv \int_D(\vert\nabla v\vert^2+\tau^2\vert v\vert^2)dx
-\int_{\partial D}c(x,\tau)\vert v\vert^2dS_x.
$$
Since $\overline B\cap \overline D=\emptyset$, $f(x)=0$ a.e. $x\in D$
and thus (2.13) gives another expression:
$$\displaystyle
J(\tau)
=\int_{\partial D}\left(\frac{\partial v}{\partial\nu}-c(x,\tau)v\right)vdS.
\tag {2.35}
$$

\proclaim{\noindent Lemma 2.2.}
There exists a positive number $\mu$ such that
$$\displaystyle
\lim_{\tau\longrightarrow\infty}\tau^{\mu}e^{2\tau\text{dist}\,(D,B)}J(\tau)>0.
\tag {2.36}
$$
\endproclaim

{\it\noindent Proof.}
Since
$$\displaystyle
\frac{\partial v}{\partial\nu}(x)
=-\frac{\tau}{4\pi}
\int_B\frac{(x-y)\cdot\nu(x)}{\vert x-y\vert^2}
e^{-\tau\vert x-y\vert}f(y)dy
-\frac{1}{4\pi}
\int_B\frac{(x-y)\cdot\nu(x)}{\vert x-y\vert^3}e^{-\tau\vert x-y\vert}f(y)dy,
$$
we have
$$
\displaystyle
J(\tau)
=\frac{\tau}{(4\pi)^2}\int_{\partial D}dS_x
\int_{B\times B}
\left(k(x,y)-\frac{1}{\tau}l(x,y)\right)
\frac{e^{-\tau(\vert x-y\vert+\vert x-y'\vert)}}{\vert x-y\vert\vert x-y'\vert}f(y)f(y')dydy'
$$
where
$$\begin{array}{c}
\displaystyle
k(x,y)=\frac{(y-x)\cdot\nu(x)}{\vert x-y\vert}-\gamma(x),
\\
\\
\displaystyle
l(x,y)=\frac{(y-x)\cdot\nu(x)}{\vert x-y\vert^2}-\beta(x).
\end{array}
$$

We divide the integrand of $J(\tau)$ into two parts.
Set $d=\text{dist}\,(\partial D, \overline{B})$
and ${\cal M}=\{(x,y)\in\partial D\times\overline B\,\vert\,\vert x-y\vert=d\}$.
It is easy to see that $d=\text{dist}\,(D, B)$.

In what follows we denote by $B_{R}(z)$ the open ball centered at a point $z$ with radius $R$.
Given $\delta>0$ define
$$\displaystyle
{\cal W}_{\delta}=\cup_{(x_0,y_0)\in{\cal M}}(\partial D\cap B_{\delta}(x_0))
\times(\overline B\cap B_{\delta}(y_0))\times(\overline B\cap B_{\delta}(y_0)).
$$
The set ${\cal W}_{\delta}$ is open in $\partial D\times\overline B\times\overline B$ and
contains the set of all $(x,y,y)$ with $(x,y)\in{\cal M}$.

\noindent
Here we state the following two claims concerning the ${\cal W}_{\delta}$.
Their proofs are almost same as those of Claims 1 and 2 in \cite{IK2}.  Just rewrite the proofs
by replacing $\overline D$ and $\partial\Omega$ in the previous definitions of ${\cal M}$, ${\cal W}_{\delta}$
and $F(x,y,y')$ in \cite{IK2} with $\partial D$ and $\overline B$; $(x-y')/\vert x-y'\vert$ in the previous
definition of $F(x,y,y')$ in \cite{IK2} with $\nu(x)$.  By this reason we omit the description.

\noindent
{\bf Claim 1.}  Given $\epsilon>0$ there exists a $\delta_1>0$ such that
for all $(x,y,y')\in {\cal W}_{\delta_1}$ it holds that
$$\displaystyle
\frac{(y-x)\cdot\nu(x)}{\vert x-y\vert}
\ge 1-\epsilon,\,\,
\vert x-y\vert\le d+\epsilon,\,\,
\vert x-y'\vert\le d+\epsilon.
$$

\noindent
{\bf Claim 2.}  Given $\delta_1>0$ there exists a $\delta_2>0$ such that
if $(x,y,y')\in\partial D\times\overline B\times\overline B\setminus{\cal W}_{\delta_1}$,
then $\vert x-y\vert+\vert x-y'\vert\ge 2d+\delta_2$.

Let $C'$ be the constant in (A1).
Give $\epsilon=C'/2$ in Claim 1 and choose $\delta_1$ in Claim 1.
Next choose $\delta_2$ in Claim 2 corresponding to $\delta_1$ already chosen.

By claim 1, we have, for all $(x,y,y')\in {\cal W}_{\delta_1}$
$k(x,y)\ge 1-\epsilon-(1-C')=C'/2$, $\vert x-y\vert\le d+\epsilon$ and $\vert x-y'\vert\le d+\epsilon$;
by claim 2, we have $e^{-\tau(\vert x-y\vert+\vert x-y'\vert)}\le e^{-2\tau d}e^{-2\tau\delta_2}$
if $(x,y,y')\in\partial D\times\overline B\times\overline B\setminus{\cal W}_{\delta_1}$.
These together with estimate $f(y)f(y')\ge C^2$ for a.e. $y\in B$ gives
$$\begin{array}{c}
\displaystyle
J(\tau)
=
\frac{\tau}{(4\pi)^2}
\int_{{\cal W}_{\delta_1}}dS_xdydy'
\left(k(x,y)-\frac{1}{\tau}l(x,y)\right)
\frac{e^{-\tau(\vert x-y\vert+\vert x-y'\vert)}}{\vert x-y\vert\vert x-y'\vert}f(y)f(y')\\
\\
\displaystyle
+
\frac{\tau}{(4\pi)^2}
\int_{\partial D\times B\times B\setminus{\cal W}_{\delta_1}}dS_xdydy'
\left(k(x,y)-\frac{1}{\tau}l(x,y)\right)
\frac{e^{-\tau(\vert x-y\vert+\vert x-y'\vert)}}{\vert x-y\vert\vert x-y'\vert}f(y)f(y')\\
\\
\displaystyle
\ge \frac{(C'/2
-\tau^{-1}C_1)\tau C^2}{(4\pi)^2(d+\epsilon)^2}
\int_{{\cal W}_{\delta_1}}dS_xdydy'
e^{-\tau(\vert x-y\vert+\vert x-y'\vert)}
+O(\tau e^{-2\tau d}e^{-\tau\delta_2}),
\end{array}
\tag {2.37}
$$
where
$$\displaystyle
C_1=\frac{1}{\text{dist}\,(D,B)}+
\Vert\beta\Vert_{L^{\infty}(\partial D)}.
$$

Choose $(x_0,y_0)\in {\cal M}$.  Since $\vert x_0-y_0\vert=d$, we have
$$\begin{array}{c}
\displaystyle
\vert x-y\vert+\vert x-y'\vert
\le 2\vert x-x_0\vert+\vert x_0-y\vert+\vert x_0-y'\vert\\
\\
\displaystyle
\le 2\vert x-x_0\vert+2d+\vert y_0-y\vert+\vert y_0-y'\vert.
\end{array}
$$
From the definition of ${\cal W}_{\delta_1}$ we have
$$\begin{array}{c}
\displaystyle
\int_{{\cal W}_{\delta_1}}dS_xdydy'
e^{-\tau(\vert x-y\vert+\vert x-y'\vert)}
\ge
\int_{\partial D\cap B_{\delta_1}(x_0)}dS_x
\int_{B\cap B_{\delta_1}(y_0)}dy
\int_{B\cap B_{\delta_1}(y_0)}dy'
e^{-\tau(\vert x-y\vert+\vert x-y'\vert)}\\
\\
\displaystyle
\ge
e^{-2\tau d}\int_{\partial D\cap B_{\delta_1}(x_0)}e^{-2\tau\vert x-x_0\vert}dS_x
\left(\int_{B\cap B_{\delta_1}(y_0)}e^{-\tau\vert y_0-y\vert}dy\right)^2.
\end{array}
\tag {2.38}
$$
Here we make use of two facts which are essentially same as Claims 3 and 4 in \cite{IK2}:

for all $\delta>0$ we have
$$\displaystyle
\liminf_{\tau\longrightarrow\infty}\tau^3\int_{B\cap B_{\delta}(y_0)}
e^{-\tau\vert y_0-y\vert}dy>0
$$
and
$$\displaystyle
\liminf_{\tau\longrightarrow\infty}\tau^2\int_{\partial D\cap B_{\delta}(x_0)}
e^{-2\tau\vert x-x_0\vert}dS_x>0.
$$
From these, (2.37) and (2.38) we see that (2.36) for $\mu=7$ is true.

\noindent
$\Box$

It follows from (2.26) that
$$\begin{array}{c}
\displaystyle
\int_{\Bbb R^3\setminus\overline D}\vert\nabla\epsilon\vert^2dx
+\tau^2\int_{\Bbb R^3\setminus\overline D}\vert\epsilon\vert^2dx+\int_{\partial D}\beta\vert\epsilon\vert^2dS
\\
\\
\displaystyle
\ge
(1-K\Vert\beta\Vert_{L^{\infty}(\partial D)}\eta^2)
\int_{\Omega\setminus\overline D}\vert\nabla\epsilon\vert^2dx+\int_{\Bbb R^3\setminus\overline\Omega}\vert\nabla\epsilon\vert^2dx\\
\\
\displaystyle
+(\tau^2-K\Vert\beta\Vert_{L^{\infty}(\partial D)}\eta^{-2})\int_{\Omega\setminus\overline D}\vert\epsilon\vert^2dx
+\tau^2\int_{\Bbb R^3\setminus\overline\Omega}\vert\epsilon\vert^2dx\\
\\
\displaystyle
\ge
(\tau^2-K\Vert\beta\Vert_{L^{\infty}(\partial D)}\eta^{-2})\int_{\Bbb R^3\setminus\overline D}\vert\epsilon\vert^2dx
\\
\\
\displaystyle
\ge
\frac{\tau^2}{2}\int_{\Bbb R^3\setminus\overline D}\vert\epsilon\vert^2dx
\end{array}
\tag {2.39}
$$
for a fixed $\eta$ with $1\ge K\Vert\beta\Vert_{L^{\infty}(\partial D)}\eta^2$ and $\tau\ge \sqrt{2K\Vert\beta\Vert_{L^{\infty}(\partial D)}}\eta^{-1}$.
Moreover we have
$$\begin{array}{c}
\displaystyle
\frac{\tau^2}{2}\int_{\Bbb R^3\setminus\overline D}\vert\epsilon\vert^2dx
+e^{-\tau T}\int_{\Bbb R^3\setminus\overline D}F\epsilon dx
=\frac{1}{2}\int_{\Bbb R^3\setminus\overline D}
\left\vert\tau\epsilon+\frac{e^{-\tau T}}{\tau}F\right\vert^2dx
-\frac{e^{-2\tau T}}{2\tau^2}
\int_{\Bbb R^3\setminus\overline D}\vert F\vert^2 dx\\
\\
\displaystyle
\ge
-\frac{e^{-2\tau T}}{2\tau^2}
\int_{\Bbb R^3\setminus\overline D}\vert F\vert^2 dx
\end{array}
\tag {2.40}
$$
and
$$\begin{array}{c}
\displaystyle
\tau\int_{\partial D}\gamma\vert\epsilon\vert^2dS
+e^{-\tau T}\int_{\partial D}G\epsilon dS\\
\\
\displaystyle
=\tau\int_{\partial D}\gamma\left\vert\epsilon+\frac{e^{-\tau T}}{2\tau}u(x,T)\right\vert^2dS
-\frac{e^{-2\tau T}}{4\tau}\int_{\partial D}\gamma\vert u(x,T)\vert^2dx\\
\\
\displaystyle
\ge
-\frac{e^{-2\tau T}}{4\tau}\int_{\partial D}\gamma\vert u(x,T)\vert^2dx.
\end{array}
\tag {2.41}
$$
Note that we have made use of the assumption $\gamma\ge 0$ on $\partial D$.
Now it follows from (2.15), (2.39), (2.40) and (2.41) that
$$\begin{array}{c}
\displaystyle
I(\tau)
=J(\tau)
+\int_{\Bbb R^3\setminus\overline D}\vert\nabla\epsilon\vert^2dx+\tau^2\int_{\Bbb R^3\setminus\overline D}\vert\epsilon\vert^2dx
+\int_{\partial D}c(x,\tau)\vert\epsilon\vert^2dS\\
\\
\displaystyle
+e^{-\tau T}\int_{\Bbb R^3\setminus\overline D}F\epsilon dx-e^{-\tau T}\int_{\Omega\setminus\overline D}Fvdx
+e^{-\tau T}\int_{\partial D}G\epsilon dS\\
\\
\displaystyle
\ge J(\tau)
-\frac{e^{-2\tau T}}{2\tau^2}
\int_{\Bbb R^3\setminus\overline D}\vert F\vert^2 dx
-\frac{e^{-2\tau T}}{4\tau}\int_{\partial D}\gamma\vert u(x,T)\vert^2dx
-e^{-\tau T}\int_{\Omega\setminus\overline D}Fvdx\\
\\
\displaystyle
=J(\tau)+O(e^{-2\tau T}+\tau e^{-\tau(T+\text{dist}\,(\Omega,\,B))}).
\end{array}
$$
and thus
$$\begin{array}{c}
\displaystyle
\tau^{\mu}e^{2\tau\text{dist}\,(D,\,B)}I(\tau) \ge
\tau^{\mu}e^{2\tau\text{dist}\,(D,\,B)}J(\tau)\\
\\
\displaystyle
+O(\tau^{\mu}e^{-2\tau(T-\text{dist}\,(D,\,B))}
+\tau^{\mu+1}
e^{-\tau(T-2\text{dist}\,(D,\,B)+\text{dist}\,(\Omega,\,B))}).
\end{array}
$$
Now from this together with (2.36) we obtain (2.34).

\subsection{Case (A2).  An estimate of $I(\tau)$ from above}

In this subsection we prove that there exists a positive constant $\mu$ such that
$$\displaystyle
\liminf_{\tau\longrightarrow\infty}\left(-\tau^{\mu}e^{2\tau\text{dist}\,(D,\,B)}I(\tau)\right)>0
\tag {2.42}
$$
provided $\gamma$ satisfies (A2).

Having (2.16) and (2.42), we immediately obtain (1.6)
and (1.7).

The proof of (2.42) starts with rewriting (2.12) as
$$\begin{array}{c}
\displaystyle
\int_{\partial D}\frac{1}{4c(x,\tau)}\left\vert\left(\frac{\partial v}{\partial\nu}-c(x,\tau)v\right)-e^{-\tau T}G\right\vert^2dS
+\frac{e^{-2\tau T}}{4\tau^2}\int_{\Bbb R^3\setminus\overline D}
\vert F\vert^2dx\\
\\
\displaystyle
=\int_{\partial D}c(x,\tau)\left\vert\epsilon-\frac{1}{2c(x,\tau)}
\left\{\left(\frac{\partial v}{\partial\nu}-c(x,\tau)v\right)-e^{-\tau T}G\right\}\right\vert^2dS\\
\\
\displaystyle
+\int_{\Bbb R^3\setminus\overline D}
\left(\vert\nabla\epsilon\vert^2
+\tau^2\left\vert
\epsilon+\frac{e^{-\tau T}}{2\tau^2}F\right\vert^2\right)dx.
\end{array}
\tag {2.43}
$$
Since
$$\displaystyle
\left\vert
\epsilon+\frac{e^{-\tau T}}{2\tau^2}F\right\vert^2
\ge\frac{1}{2}\vert\epsilon\vert^2-
\frac{e^{-2\tau T}}{4\tau^4}\vert F\vert^2
$$
and
$$\begin{array}{c}
\displaystyle
\left\vert\epsilon-\frac{1}{2c(x,\tau)}
\left\{\left(\frac{\partial v}{\partial\nu}-c(x,\tau)v\right)-e^{-\tau T}G\right\}\right\vert^2\\
\\
\displaystyle
\ge
\frac{1}{2}\vert\epsilon\vert^2
-\frac{1}{4c(x,\tau)^2}\left\vert\left(\frac{\partial v}{\partial\nu}-c(x,\tau)v\right)-e^{-\tau T}G\right\vert^2,
\end{array}
$$
it follows from (2.43) that
$$\begin{array}{c}
\displaystyle
\int_{\partial D}\frac{c(x,\tau)}{2}\vert\epsilon\vert^2dS
+\int_{\Bbb R^3\setminus\overline D}
\left(\vert\nabla\epsilon\vert^2+\frac{\tau^2}{2}\vert\epsilon\vert^2\right)dx\\
\\
\displaystyle
\le
\int_{\partial D}\frac{1}{2c(x,\tau)}\left\vert\left(\frac{\partial v}{\partial\nu}-c(x,\tau)v\right)-e^{-\tau T}G\right\vert^2dS
+\frac{e^{-2\tau T}}{2\tau^2}\int_{\Bbb R^3\setminus\overline D}
\vert F\vert^2dx
\end{array}
$$
and thus
$$\begin{array}{c}
\displaystyle
\int_{\partial D}c(x,\tau)\vert\epsilon\vert^2dS
+\int_{\Bbb R^3\setminus\overline D}
\left(\vert\nabla\epsilon\vert^2+\tau^2\vert\epsilon\vert^2\right)dx\\
\\
\displaystyle
\le
\int_{\partial D}\frac{1}{c(x,\tau)}\left\vert\left(\frac{\partial v}{\partial\nu}-c(x,\tau)v\right)-e^{-\tau T}G\right\vert^2dS
+\frac{e^{-2\tau T}}{\tau^2}\int_{\Bbb R^3\setminus\overline D}
\vert F\vert^2dx\\
\\
\displaystyle
=\int_{\partial D}\frac{1}{c(x,\tau)}\left\vert\frac{\partial v}{\partial\nu}-c(x,\tau)v\right\vert^2dS
+O(e^{-\tau(T+\text{dist}\,(D,\,B))}+e^{-2\tau T}).
\end{array}
\tag {2.44}
$$
It follows from (2.17) that
$$\displaystyle
e^{-\tau T}\int_{\Bbb R^3\setminus\overline D}F\epsilon dx=
O(\tau(e^{-2\tau T}+e^{-\tau(T+\text{dist}\,(D,\,B))}))
\tag {2.45}
$$
and from (2.28) that
$$\displaystyle
e^{-\tau T}\int_{\partial D}G\epsilon dS
=O(\tau^{1/2}(e^{-2\tau T}+e^{-\tau(T+\text{dist}\,(D,\,B))}).
\tag {2.46}
$$
Now applying (2.31), (2.32), (2.44), (2.45) and (2.46) to (2.15) we obtain
$$\begin{array}{c}
\displaystyle
I(\tau)\le
\tilde{J}(\tau)
+O(\tau(e^{-2\tau T}+e^{-\tau(T+\text{dist}\,(D,\,B))}))
+O(\tau e^{-\tau(T+\text{dist}\,(\Omega,\,B))}),
\end{array}
\tag {2.47}
$$
where
$$\displaystyle
\tilde{J}(\tau)
=J(\tau)
+\int_{\partial D}\frac{1}{c(x,\tau)}\left\vert\frac{\partial v}{\partial\nu}-c(x,\tau)v\right\vert^2dS.
$$
From (2.35) one gets
$$
\begin{array}{c}
\displaystyle
\tilde{J}(\tau)
=\int_{\partial D}\frac{1}{c(x,\tau)}\left(\frac{\partial v}{\partial\nu}-c(x,\tau)v\right)\frac{\partial v}{\partial\nu}dS\\
\\
\displaystyle
=\left(\frac{\tau}{4\pi}\right)^2
\int_{\partial D}\frac{1}{c(x,\tau)}dS_x\\
\\
\displaystyle
\times
\int_{B\times B}
\left(k(x,y')-\frac{1}{\tau}l(x,y')\right)
\left(\tilde{k}(x,y)-\frac{1}{\tau}\tilde{l}(x,y)\right)
\frac{e^{-\tau(\vert x-y'\vert+\vert x-y\vert)}}{\vert x-y'\vert\vert x-y\vert}f(y')f(y)dy'dy,
\end{array}
$$
where $k$ is the same as that of Lemma 2.2 and
$$\begin{array}{c}
\displaystyle
\tilde{k}(x,y)=\frac{(y-x)\cdot\nu(x)}{\vert x-y\vert},
\\
\\
\displaystyle
\tilde{l}(x,y)=\frac{(y-x)\cdot\nu(x)}{\vert x-y\vert^2}.
\end{array}
$$
From (A2) we have $k(x,y')\le -C'$ and thus
$$\displaystyle
k(x,y')-\frac{1}{\tau}l(x,y')\le-\frac{C'}{2},\,\,(x,y')\in\partial D\times B
$$
for all $\tau\ge\tau_0$ with $\tau_0>>1$.
Now give $\epsilon=1/2$ in Claim I of Lemma 2.2 and choose $\delta_1$ and corresponding $\delta_2$ in Claim 2.
It is easy to see that one can apply the same argument as done in Lemma 2.2 to $\tilde{J}(\tau)$ where
$\tilde{k}(x,y)-(1/\tau)\tilde{l}(x,y)$ plays the same role
as $k(x,y)-(1/\tau)l(x,y)$ in Lemma 2.2.
Thus we conclude that
there exists $\mu>0$ such that
$$\displaystyle
\liminf_{\tau\longrightarrow\infty}\left\{-\tau^{\mu}e^{2\tau\text{dist}\,(D,\,B)}\tilde{J}(\tau)\right\}>0
$$
and thus from (2.47) one gets (2.42).

\subsection{Proof of Theorem 1.2}

Since $\overline B\cap\overline D=\emptyset$, one can find a bounded open set
$\Omega$ with a smooth boundary such that $\Bbb R^3\setminus\overline\Omega$ is connected,
$\overline D\subset\Omega$ and $\overline B\cap\overline\Omega=\emptyset$.
We consider $I(\tau)$ for this $\Omega$.

Since $\text{supp}\,f\cap\overline\Omega=\emptyset$, it follows from (2.3) that
$$\displaystyle
I(\tau)
=\int_{\partial D}\frac{\partial v}{\partial\nu}wdS-\int_{\partial D}c(x,\tau)wvdS
-e^{-\tau T}\int_{\partial D}GvdS-e^{-\tau T}\int_{\Omega\setminus\overline D}Fvdx.
\tag {2.48}
$$
Substituting $\phi=w$ into (2.10), we have
$$
\displaystyle
\int_{\partial D}\frac{\partial v}{\partial\nu}w dS
=\int_{\Bbb R^3\setminus\overline D}fw dx-\int_{\Bbb R^3\setminus\overline D}(\nabla v\cdot\nabla w
+\tau^2vw)dx
$$
and thus (2.48) becomes
$$\begin{array}{c}
\displaystyle
I(\tau)
=\int_{\Bbb R^3\setminus\overline D}f(w-v)dx-\left\{\int_{\Bbb R^3\setminus\overline D}(\nabla w\cdot\nabla v
+(\tau^2w-f)v)dx+\int_{\partial D}c(x,\tau)wvdS\right\}\\
\\
\displaystyle
-e^{-\tau T}\int_{\partial D}GvdS-e^{-\tau T}\int_{\Omega\setminus\overline D}Fvdx.
\end{array}
$$
Now it follows from this and (2.2) for $\phi=v\vert_{\Bbb R^3\setminus\overline D}$ that
$$\displaystyle
I(\tau)
=\int_{\Bbb R^3\setminus\overline D}f(w-v)dx
+e^{-\tau T}\int_{\Bbb R^3\setminus\overline\Omega}Fvdx.
\tag {2.49}
$$
Since $\Vert F\Vert_{L^2(\Bbb R^3\setminus\overline D)}=O(\tau)$ and
$$\displaystyle
\Vert v\Vert_{L^2(\Bbb R^3)}\le\tau^{-2}\Vert f\Vert_{L^2(\Bbb R^3)},
\tag {2.50}
$$
(2.49) yields
$$\displaystyle
\int_{\Bbb R^3\setminus\overline D}f(w-v)dx=I(\tau)+O(\tau^{-1}e^{-\tau T}).
\tag {2.51}
$$
Now (1.10) and (1.12) follow from (2.51) together with (2.16) and (2.34) in case (A1);
(1.11) and (1.12) follow from (2.51) together with (2.16) and (2.42) in case (A2).

\section{Proof of Theorems 1.3 and 1.4}

Set $V=H^1(\Bbb R^3)$ and $H=L^2(\Bbb R^3)$. Applying Theorem 1 on
p 558 in \cite{DL} to (1.13), we know that: given $u^0\in V$ and
$u^1\in H$ there exists a unique $u$ satisfying
$$
\displaystyle
u\in L^2(0,\,\,T;V),\,\,
u'\in L^2(0,\,\,T;V),\,\,
d/dt(C(u'(\,\cdot\,))\in L^2(0,\,\,T;V' )
$$
such that, for all $\phi\in V$
$$\displaystyle
<\frac{d}{dt}C(u'(t)),\phi>
+a(u(t),\phi)=0,\,\,\text{a.e.}\,t\in]0,\,\,T[,
\tag {3.1}
$$
and $u(0)=u^0$ and $u'(0)=u^1$, where
$$\displaystyle
a(u,v)=\int_{\Bbb R^3}\nabla u\cdot\nabla vdx,\,\,u,v\in V
$$
and $C:H\longrightarrow H$ is the bounded linear operator defined
by
$$\displaystyle
C(u)=\alpha u,\,\,u\in H.
$$
Note that this $C$ satisfies (5.11) on p 553 in \cite{DL}
under the condition $\alpha(x)\ge C$ a.e. $x\in\Bbb R^3$
for some positive constant.  Note that at this step we do not make use
of condition $\alpha(x)=1$ a.e. $x\in\Bbb R^3\setminus D$.

Note also that this is the case when $b=0$ in Theorem 1. However,
the equation is {\it homogeneous}, that is, $f=0$ in their
notation and by virtue of this their proof covers also this case.

In this section we say that this $u$ for $u^0=0$ and
$u^1=f$ is the weak solution of (1.13).

We see that $w$ given by (1.14) belongs to $V$ and
applying integration by parts to (3.1) multiplied by $e^{-\tau T}$,
we obtain, for all $\phi\in V$
$$\displaystyle
\int_{\Bbb R^3}(\nabla w\cdot\nabla\phi+\tau^2\alpha w\phi)dx
-\int_{\Bbb R^3}\alpha f\phi dx+\int_{\Bbb R^3}\mbox{\boldmath $F$}\phi dx=0,
\tag {3.2}
$$
where
$$\displaystyle
\mbox{\boldmath $F$}(x)=e^{-\tau T}\alpha(x)(u'(x,T)+\tau u(x,T)).
$$
This means that $w$ is the weak solution of the following equation
$$\displaystyle
(\triangle-\alpha\tau^2)w+\alpha f
=\mbox{\boldmath $F$}\,\,\text{in}\,\Bbb R^3.
\tag {3.3}
$$
By the same reason as $w$ in Section 2 we have $w\in H^2_{\text{loc}}(\Bbb R^3)$
and (3.3) holds a.e. $x\in\Bbb R^3$.

Since $v$ satisfies (1.3) in the weak sense,
we see that $w-v$ satisfies, in a weak sense,
$$\displaystyle
(\triangle-\alpha\tau^2)(w-v)
=(\alpha-1)(\tau^2 v-f)
+\mbox{\boldmath $F$}\,\,\text{in}\,\Bbb R^3.
\tag {3.4}
$$

\subsection{Two expressions}

We start with the following two key expressions.
Note that in the proposition we do not assume that
$\text{supp}\,f\cap\overline\Omega=\emptyset$.

\proclaim{\noindent Proposition 3.1.}
We have:
$$\begin{array}{c}
\displaystyle
\int_{\partial\Omega}\left(\frac{\partial v}{\partial\nu}w
-\frac{\partial w}{\partial\nu}v\right)dS\\
\\
\displaystyle
=\int_{\Bbb R^3}\left\{\vert\nabla(w-v)\vert^2+\tau^2\alpha\vert w-v\vert^2\right\}dx
+\tau^2\int_{D}(1-\alpha)\vert v\vert^2dx\\
\\
\displaystyle
+\int_{\Bbb R^3}(1-\alpha)f(w-v)dx
+\int_{\Bbb R^3}\mbox{\boldmath $F$}(w-v)dx
+\int_{\Omega}f(\alpha v-w)dx-\int_{\Omega}\mbox{\boldmath $F$}vdx;
\end{array}
\tag {3.5}
$$
$$\begin{array}{c}
\displaystyle
\int_{\partial\Omega}\left(\frac{\partial w}{\partial\nu}v-\frac{\partial v}{\partial\nu}w\right)dS\\
\\
\displaystyle
=\int_{\Bbb R^3}\left\{\vert\nabla(v-w)\vert^2+\tau^2\vert v-w\vert^2\right\}dx
+\tau^2\int_{\Omega}(\alpha-1)\vert w\vert^2 dx\\
\\
\displaystyle
+\int_{\Bbb R^3}(\alpha-1)f(v-w)dx
-\int_{\Bbb R^3}\mbox{\boldmath $F$}(v-w)dx
+\int_{\Omega}f(w-\alpha v)dx
+\int_{\Omega}\mbox{\boldmath $F$}vdx.
\end{array}
\tag {3.6}
$$

\endproclaim

{\it\noindent Proof.}
The proof of (3.5) and (3.6) are rather simpler than that of (2.3) and (2.4).
First we give a proof of (3.5).
Since $v\in H^2(\Omega)$ and satisfies $(\triangle-\tau^2)v+f=0$ a.e. $x\in\Omega$,
integration by parts gives
$$
\displaystyle
\int_{\partial\Omega}\frac{\partial v}{\partial\nu}w
dS
=\int_{\Omega}(\tau^2 v-f)wdx+\int_{\Omega}\nabla v \cdot\nabla w dx
$$
and by the same reason for $w$ we have
$$\displaystyle
\int_{\partial\Omega}\frac{\partial w}{\partial\nu}vdS
=\int_{\Omega}\nabla w\cdot\nabla vdx+\tau^2\int_{\Omega}\alpha w vdx
-\int_{\Omega}\alpha f vdx+\int_{\Omega}\mbox{\boldmath $F$}vdx.
$$
These yield
$$\begin{array}{c}
\displaystyle
\int_{\partial\Omega}\left(\frac{\partial v}{\partial\nu}w-\frac{\partial w}{\partial\nu}v\right)dS
=\tau^2\int_{\Omega}(1-\alpha)wvdx+\int_{\Omega}f(\alpha v-w)dx
-\int_{\Omega}\mbox{\boldmath $F$}vdx.
\end{array}
\tag {3.7}
$$
Write the integral in the first term of
the right hand-side as
$$\displaystyle
\int_{\Omega}(1-\alpha)wvdx=\int_{\Omega}(1-\alpha)v(w-v)dx+\int_{\Omega}(1-\alpha)\vert v\vert^2 dx.
\tag {3.8}
$$
On the other hand, from the weak form of (3.4) and integration by parts we have
$$\begin{array}{c}
\displaystyle
-\int_{\Bbb R^3}\vert\nabla(w-v)\vert^2dx-\tau^2\int_{\Bbb R^3}\alpha\vert w-v\vert^2dx\\
\\
\displaystyle
=\tau^2\int_{\Bbb R^3}(\alpha-1)v(w-v)dx-\int_{\Bbb R^3}(\alpha-1)f(w-v)dx
+\int_{\Bbb R^3}\mbox{\boldmath $F$}(w-v)dx,
\end{array}
$$
that is,
$$\begin{array}{c}
\displaystyle
\tau^2\int_{\Bbb R^3}(1-\alpha)v(w-v)dx\\
\\
\displaystyle
=
\int_{\Bbb R^3}\{\vert\nabla(w-v)\vert^2+\tau^2\alpha\vert w-v\vert^2\}dx
+\int_{\Bbb R^3}(1-\alpha)f(w-v)dx
+\int_{\Bbb R^3}\mbox{\boldmath $F$}(w-v)dx.
\end{array}
\tag {3.9}
$$
Since $1-\alpha(x)=0$ outside $D$, from (3.7), (3.8) and (3.9) we obtain the desired expression.

Second we give a proof of (3.6).

We change the role of $v$ and $w$.
Introducing
$$\displaystyle
\tilde{f}=f-\frac{1}{\alpha}\mbox{\boldmath $F$},\,\,
\tilde{\mbox{\boldmath $F$}}=-\frac{1}{\alpha}\mbox{\boldmath $F$},
$$
one can rewrite equations (1.4) and (3.3) as
$$\displaystyle
(\triangle-\tau^2)v+\tilde{f}=\tilde{\mbox{\boldmath $F$}}\,\,\text{in}\,\Bbb R^3
$$
and
$$\displaystyle
(\triangle-\alpha\tau^2)w+\alpha\tilde{f}=0\,\,\text{in}\,\Bbb R^3.
$$
Then $v-w$ satisfies
$$\begin{array}{c}
\displaystyle
(\triangle-\tau^2)(v-w)
=\tilde{\mbox{\boldmath $F$}}-\tilde{f}-(\triangle-\tau^2)w\\
\\
\displaystyle
=\tilde{\mbox{\boldmath $F$}}-\tilde{f}-\{\triangle-\alpha\tau^2+(\alpha-1)\tau^2\}w\\
\\
\displaystyle
=\tilde{\mbox{\boldmath $F$}}-\tilde{f}+\alpha\tilde{f}+(1-\alpha)\tau^2w\\
\\
\displaystyle
=(1-\alpha)(\tau^2w-\tilde{f})+\tilde{\mbox{\boldmath $F$}}.
\end{array}
\tag {3.10}
$$
Integration by parts gives
$$\begin{array}{c}
\displaystyle
\int_{\partial\Omega}\frac{\partial w}{\partial\nu}vdS
=\int_{\Omega}(\triangle w)vdx+
\int_{\Omega}\nabla w\cdot\nabla v dx\\
\\
\displaystyle
=\int_{\Omega}\alpha\tau^2 wvdx-\int_{\Omega}\alpha\tilde{f}vdx
+\int_{\Omega}\nabla w\cdot\nabla vdx\\
\\
\displaystyle
=\int_{\Omega}\alpha(\tau^2w-\tilde{f})vdx+\int_{\Omega}\nabla w\cdot\nabla vdx;
\end{array}
$$
and
$$\begin{array}{c}
\displaystyle
\int_{\partial\Omega}\frac{\partial v}{\partial\nu}wdS
=\int_{\Omega}(\triangle v)wdx+\int_{\Omega}\nabla v\cdot\nabla w dx\\
\\
\displaystyle
=\int_{\Omega}(\tau^2v-\tilde{f}+\tilde{\mbox{\boldmath $F$}})wdx
+\int_{\Omega}\nabla v\cdot\nabla w dx\\
\\
\displaystyle
=\int_{\Omega}\nabla v\cdot\nabla w dx+\tau^2\int_{\Omega} vwdx
-\int_{\Omega}\tilde{f}wdx+\int_{\Omega}\tilde{\mbox{\boldmath $F$}}w dx.
\end{array}
$$
From these we obtain
$$\begin{array}{c}
\displaystyle
\int_{\partial\Omega}\left(\frac{\partial w}{\partial\nu}v-\frac{\partial v}{\partial\nu}w\right)dS\\
\\
\displaystyle
=\tau^2\int_{\Omega}(\alpha-1)wvdx
+\int_{\Omega}\tilde{f}(w-\alpha v)dx
-\int_{\Omega}\tilde{\mbox{\boldmath $F$}}wdx.
\end{array}
\tag {3.11}
$$
Here rewrite the integral of the first term in the right-hand side:
$$\displaystyle
\int_{\Omega}(\alpha-1)wvdx
=\int_{\Omega}(\alpha-1)w(v-w)dx
+\int_{\Omega}(\alpha-1)\vert w\vert^2 dx.
\tag {3.12}
$$
It follows from the weak form of (3.10) that
$$\begin{array}{c}
\displaystyle
-\int_{\Bbb R^3}\vert\nabla(v-w)\vert^2dx-\tau^2\int_{\Bbb R^3}\vert v-w\vert^2 dx\\
\\
\displaystyle
=\int_{\Bbb R^3}(1-\alpha)(\tau^2w-\tilde{f})(v-w)dx
+\int_{\Bbb R^3}\tilde{\mbox{\boldmath $F$}}(v-w)dx.
\end{array}
$$
This yields
$$\begin{array}{c}
\displaystyle
\tau^2\int_{\Bbb R^3}(\alpha-1)w(v-w)dx\\
\\
\displaystyle
=\int_{\Bbb R^3}\{\vert\nabla(v-w)\vert^2+\tau^2\vert v-w\vert^2\}dx
-\int_{\Bbb R^3}(1-\alpha)\tilde{f}(v-w)dx
+\int_{\Bbb R^3}\tilde{\mbox{\boldmath $F$}}(v-w)dx\\
\\
\displaystyle
=\int_{\Bbb R^3}\{\vert\nabla(v-w)\vert^2+\tau^2\vert v-w\vert^2\}dx
+\int_{\Bbb R^3}\{\tilde{\mbox{\boldmath $F$}}+(\alpha-1)\tilde{f}\}(v-w)dx.
\end{array}
\tag {3.13}
$$
Now from (3.11), (3.12) and (3.13) we obtain
$$\begin{array}{c}
\displaystyle
\int_{\partial\Omega}\left(\frac{\partial w}{\partial\nu}v-\frac{\partial v}{\partial\nu}w\right)dS
\\
\\
\displaystyle
=\tau^2\int_{\Omega}(\alpha-1)\vert w\vert^2 dx
+\int_{\Bbb R^3}\{\vert\nabla(v-w)\vert^2+\tau^2\vert v-w\vert^2\}dx\\
\\
\displaystyle
+\int_{\Bbb R^3}\{\tilde{\mbox{\boldmath $F$}}+(\alpha-1)\tilde{f}\}(v-w)dx
+\int_{\Omega}\tilde{f}(w-\alpha v)dx-\int_{\Omega}\tilde{\mbox{\boldmath $F$}}wdx.
\end{array}
$$
Since
$$\begin{array}{c}
\displaystyle
\tilde{\mbox{\boldmath $F$}}+(\alpha-1)\tilde{f}=-\mbox{\boldmath $F$}+(\alpha-1)f,\\
\\
\displaystyle
\tilde{f}-\tilde{\mbox{\boldmath $F$}}=f,
\end{array}
$$
this is nothing but (3.6).

\noindent
$\Box$

Now assume that $\text{supp}\,f\cap\overline\Omega=\emptyset$.
Set
$$\displaystyle
K(\tau)=\int_{\partial\Omega}\left(\frac{\partial v}{\partial\nu}w
-\frac{\partial w}{\partial\nu}v\right)dS,\,\,\tau>0.
$$
Since $\alpha(x)=1$ a. e. $x\in\Bbb R^3\setminus D$, (3.5) and
(3.6) become the following two expressions of $K(\tau)$:
$$\begin{array}{c}
\displaystyle
K(\tau)
=\int_{\Bbb R^3}\left\{\vert\nabla(w-v)\vert^2+\tau^2\alpha\vert w-v\vert^2\right\}dx\\
\\
\displaystyle
+\tau^2\int_{D}(1-\alpha)\vert v\vert^2dx+\int_{\Bbb R^3}\mbox{\boldmath $F$}(w-v)dx
-\int_{\Omega}\mbox{\boldmath $F$}vdx;
\end{array}
\tag {3.14}
$$
$$\begin{array}{c}
\displaystyle
-K(\tau)
=\int_{\Bbb R^3}\left\{\vert\nabla(v-w)\vert^2+\tau^2\vert v-w\vert^2\right\}dx
+\tau^2\int_{\Omega}(\alpha-1)\vert w\vert^2 dx\\
\\
\displaystyle
-\int_{\Bbb R^3}\mbox{\boldmath $F$}(v-w)dx
+\int_{\Omega}\mbox{\boldmath $F$}vdx.
\end{array}
\tag {3.15}
$$

\subsection{An estimate of $\vert K(\tau)\vert$ from above}

We describe the proof of the following estimate as
$\tau\longrightarrow\infty$:
$$
\displaystyle
e^{2\tau\text{dist}\,(D,\,B)}\vert K(\tau)\vert
=O(\tau^2).
\tag {3.16}
$$

Since $(\alpha(x)-1)f(x)=0$ for a.e. $x\in\Bbb R^3$ it follows from (3.9) that
$$\begin{array}{c}
\displaystyle
\int_{\Bbb R^3}\left\{\vert\nabla(w-v)\vert^2+\tau^2\alpha\vert w-v\vert^2\right\}dx\\
\\
\displaystyle
=\tau^2\int_{D}(1-\alpha)v(w-v)dx-\int_{\Bbb R^3}\mbox{\boldmath $F$}(w-v)dx.
\end{array}
\tag {3.17}
$$
Applying a much simpler argument than that for the derivation of (2.17) and (2.18) in Lemma 2.1
to (3.17) and using the estimate
$$\displaystyle
\Vert\mbox{\boldmath $F$}\Vert_{L^2(\Bbb R^3)}^2
=O(\tau^2 e^{-2\tau T}),
\tag {3.18}
$$
we obtain, as $\tau\longrightarrow\infty$
$$\displaystyle
\begin{array}{c}
\Vert w-v\Vert_{L^2(\Bbb R^3)}^2=O(e^{-2\tau T}+e^{-2\tau\text{dist}\,(D,\,B)}),\\
\\
\displaystyle
\Vert\nabla(w-v)\Vert_{L^2(\Bbb R^3)}^2=O(\tau^2(e^{-2\tau T}+e^{-2\tau\text{dist}\,(D,\,B)}))
\end{array}
\tag {3.19}
$$
and thus one gets
$$
\displaystyle
\left\vert\int_{\Bbb R^3}\mbox{\boldmath $F$}(w-v)dx\right\vert
=O(\tau(e^{-2\tau T}+e^{-\tau(T+\text{dist}\,(D,\,B))})).
\tag {3.20}
$$
Since $\Vert v\Vert_{L^2(\Omega)}=O(e^{-\tau\text{dist}\,(\Omega,\,B)})$, from (3.18)
one gets also
$$
\displaystyle
\left\vert\int_{\Omega}\mbox{\boldmath $F$}vdx\right\vert
=O(\tau e^{-\tau(T+\text{dist}\,(\Omega,\,B))}).
\tag {3.21}
$$
Now it follows from (3.14), $\Vert v\Vert_{L^2(D)}
=O(e^{-\tau\text{dist}\,(D,\,B)})$ and (3.18) to (3.21) that, as $\tau\longrightarrow\infty$
$$\displaystyle
\vert K(\tau)\vert
=O(\tau^2(e^{-2\tau T}+e^{-2\tau\text{dist}\,(D,\,B)}))+\tau
e^{-\tau(T+\text{dist}\,(\Omega,\,B))}).
$$
From this we obtain (3.16).

\subsection{Case (B1).  An estimate of $K(\tau)$ from below}

In this subsection we consider the case when $\alpha(x)\le 1-C'$ a.e. $x\in D$
and prove that:
$$\displaystyle
\liminf_{\tau\longrightarrow\infty}\tau^4e^{2\tau\text{dist}\,(D,\,B)}
K(\tau)>0.
\tag {3.22}
$$

Using the completing square
$$\displaystyle
\tau^2\alpha\vert w-v\vert^2
+\mbox{\boldmath $F$}(w-v)
=\tau^2\alpha\left\vert(w-v)+\frac{\mbox{\boldmath $F$}}{2\tau^2\alpha}\right\vert^2
-\frac{\vert\mbox{\boldmath $F$}\vert^2}{4\tau^2\alpha},
$$
from (3.14), (3.18) and (3.21) we obtain
$$\displaystyle
K(\tau)
\ge
C'\tau^2\int_D\vert v\vert^2dx
+O(e^{-2\tau T}+\tau e^{-\tau(T+\text{dist}\,(\Omega,\,B))}).
\tag {3.23}
$$

Here we cite the following fact established in \cite{IE0}:
$$\displaystyle
\liminf_{\tau\longrightarrow\infty}
\tau^6e^{2\tau\text{dist}\,(D,\,B)}\int_D\vert v\vert^2dx>0.
\tag {3.24}
$$
This together with (3.23) yields (3.22).

\subsection{Case (B2).  An estimate of $K(\tau)$ from above}

Next consider the case when $\alpha(x)\ge 1+C'$ a.e. $x\in D$.
In this subsection we prove the following estimate:
$$\displaystyle
\liminf_{\tau\longrightarrow\infty}\left(-\tau^4e^{2\tau\text{dist}\,(D,\,B)}K(\tau)\right)
>0.
\tag {3.25}
$$

We make use of the following completing the square:
$$\begin{array}{c}
\displaystyle
\vert v-w\vert^2+(\alpha-1)\vert w\vert^2
=\vert v\vert^2-2vw+\alpha\vert w\vert^2\\
\\
\displaystyle
=\left\vert\sqrt{\alpha}\,w-\frac{v}{\sqrt{\alpha}}\right\vert^2
+\frac{\alpha-1}{\alpha}\vert v\vert^2.
\end{array}
$$
This together with (3.15), (3.20) and (3.21) gives
$$\begin{array}{c}
\displaystyle
-K(\tau)\ge\tau^2\int_D\frac{\alpha-1}{\alpha}\vert v\vert^2dx
-\int_{\Bbb R^3}\mbox{\boldmath $F$}(v-w)dx
+\int_{\Omega}\mbox{\boldmath $F$}vdx\\
\\
\displaystyle
\ge\frac{C'\tau^2}{\Vert\alpha\Vert_{L^{\infty}(\Bbb R^3)}}
\int_D\vert v\vert^2dx
+O(\tau(e^{-2\tau T}+e^{-\tau(T+\text{dist}\,(D,\,B))}))
+O(\tau e^{-\tau(T+\text{dist}\,(\Omega,\,B))}).
\end{array}
$$
Now a combination of this and (3.24) yields (3.25).

It is easy to see that (1.15)/(1.16) and (1.17) follow
from the combination of (3.16) and (3.22)/(3.25) in the case when
(B1)/(B2) is satisfied.

\subsection{Proof of Theorem 1.4.}

Since $\overline B\cap\overline D=\emptyset$, one can find a bounded open set
$\Omega$ with a smooth boundary such that $\Bbb R^3\setminus\overline\Omega$ is connected,
$\overline D\subset\Omega$ and $\overline B\cap\overline\Omega=\emptyset$.

We make use of (3.7).  Since $\text{supp}\,f\cap\overline\Omega=\emptyset$, (3.7) becomes
$$\begin{array}{c}
\displaystyle
\int_{\partial\Omega}\left(\frac{\partial v}{\partial\nu}w-\frac{\partial w}{\partial\nu}v\right)dS
=\tau^2\int_{\Omega}(1-\alpha)wvdx
-\int_{\Omega}\mbox{\boldmath $F$}vdx.
\end{array}
\tag {3.26}
$$
One can rewrite the first term of this right-hand side as
$$\begin{array}{c}
\displaystyle
\tau^2\int_{\Omega}(1-\alpha)wvdx
=\tau^2\int_{\Bbb R^3}(1-\alpha)wvdx\\
\\
\displaystyle
=\tau^2\int_{\Bbb R^3}vwdx-\tau^2\int_{\Bbb R^3}\alpha wvdx\\
\\
\displaystyle
=-\int_{\Bbb R^3}\nabla v\cdot\nabla wdx+\int_{\Bbb R^3}fwdx
-\tau^2\int_{\Bbb R^3}\alpha wvdx\\
\\
\displaystyle
=-\int_{\Bbb R^3}(\nabla w\cdot\nabla v+\tau^2\alpha wv)dx+\int_{\Bbb R^3}fwdx\\
\\
\displaystyle
=-\left(\int_{\Bbb R^3}\alpha fvdx-\int_{\Bbb R^3}\mbox{\boldmath $F$}vdx\right)+\int_{\Bbb R^3}fwdx\\
\\
\displaystyle
=\int_{\Bbb R^3}f(w-\alpha v)dx+\int_{\Bbb R^3}\mbox{\boldmath $F$}vdx\\
\\
\displaystyle
=\int_{\Bbb R^3}f(w-v)dx+\int_{\Bbb R^3}(1-\alpha)fvdx+\int_{\Bbb R^3}\mbox{\boldmath $F$}vdx.
\end{array}
$$
Since $(1-\alpha)f=0$ in $\Bbb R^3$, from this and (3.26) we obtain
$$\displaystyle
\int_{\partial\Omega}\left(\frac{\partial v}{\partial\nu}w-\frac{\partial w}{\partial\nu}v\right)dS
=\int_{\Bbb R^3}f(w-v)dx+\int_{\Bbb R^3\setminus\overline\Omega}\mbox{\boldmath $F$}vdx.
$$
This together with (2.50) and (3.18) gives
$$\displaystyle
\int_{\Bbb R^3}f(w-v)dx
=\int_{\partial\Omega}\left(\frac{\partial v}{\partial\nu}w-\frac{\partial w}{\partial\nu}v\right)dS
+O(\tau^{-1}e^{-\tau T}).
$$

Now (1.18)/(1.19) and (1.20) follows from this
together with (3.16) and (3.22)/(3.25) in the case (B1)/(B2).

\section{Remarks and open problems}

\subsection{One space dimensional case}

Let us consider the meaning of conditions (A1) and (A2) in one
space dimensional case.

Let $0<T<\infty$ and $a>0$.
Given $f\in L^2(\Bbb R)$ with $\text{supp}\,f\subset\,]-\infty, a[$ let $u=u_f(x,t)$ be the weak solution of the following
initial boundary value problem:
$$\displaystyle\begin{array}{c}
\displaystyle
\partial_t^2u-\partial_x^2u=0\,\,\text{in}\,]-\infty,\,a[\times\,]0,\,T[,\\
\\
\displaystyle
u(x,0)=0,\,\,\partial_tu(x,0)=f(x)\,\,\text{in}\,]-\infty,\,a[,\\
\\
\displaystyle
-\partial_xu(a,t)-\gamma\partial_tu(a,t)-\beta u(a,t)=0,\,\,0<t<T,
\end{array}
\tag {4.1}
$$
where $\gamma$ and $\beta\in\Bbb R$ are constant and $\gamma\ge 0$.
This is the case when $D=]a,\,\infty[$.  We take $\Omega=]0,\,\,\infty[$ and
thus $\partial\Omega=\{0\}$.

Hereafter we choose initial data $f$ in such a way that
$\text{supp}\,f\subset\,]-\infty,\,0[$.

Define
$$\displaystyle
w=w_f(x,\tau)
=\int_0^Te^{-\tau t}u_f(x,t)dt,\,\,x\in\,]-\infty,a[,\,\,\tau>0.
$$

$w$ satisfies
$$
\displaystyle
w''-\tau^2 w+f(x)=e^{-\tau T}F(x,\tau)\,\,\text{in}\,]-\infty,\,a[
\tag {4.2}
$$
and
$$
\displaystyle
w'(a)+c(\tau)w(a)+e^{-\tau T}g(\tau)=0,
\tag {4.3}
$$
where
$$\begin{array}{c}
\displaystyle
F(x,\tau)=\partial_tu(x,T)+\tau u(x,T),\\
\\
\displaystyle
g(\tau)=\gamma u(a,T),\\
\\
\displaystyle
c(\tau)=\gamma\tau+\beta.
\end{array}
\tag {4.4}
$$

Assume that the support of $f$ is compact and satisfies
$$\displaystyle
\text{supp}\,f\subset\,[b-\epsilon,\,b]\equiv B
\tag {4.5}
$$
with $b<0$ and $\epsilon>0$.

Note that $f\in\,L^2(\Bbb R)$ with $\text{supp}\,f\subset]-\infty,\,b]$ always satisfies
$$\displaystyle
\left\vert\int_{b-\epsilon}^be^{-\tau(b-y)}f(y)dy\right\vert=O(\tau^{-1/2}).
\tag {4.6}
$$

Now define the indicator function by
$$\displaystyle
I(\tau)
=-v'(0,\tau)w(0,\tau)+w'(0,\tau)v(0,\tau),\,\,\tau>0,
$$
where
$$\displaystyle
v=v(x,\tau)
=\frac{1}{2\tau}
\int_{-\infty}^{\infty}e^{-\tau\vert x-y\vert}f(y)dy.
\tag {4.7}
$$
$v$ belongs to $H^1(\Bbb R)$ and satisfies the equation
$v''-\tau^2v+f(x)=0$ in $\Bbb R$.

We have the following asymptotic formula $I(\tau)$ as $\tau\longrightarrow\infty$.

\proclaim{\noindent Lemma 4.1.}
Let $f$ satisfy (4.5).
As $\tau\longrightarrow\infty$ we have
$$\begin{array}{c}
\displaystyle
I(\tau)
=-\frac{1}{2\tau}\frac{(\gamma-1)\tau+\beta}{(\gamma+1)\tau+\beta}
e^{-2\tau\,\text{dist}\,(D,\,B)}\left\vert\int_{b-\epsilon}^{b}e^{-\tau(b-y)}f(y)dy\right\vert^2\\
\\
\displaystyle
+O(e^{-\tau(T+\text{dist}\,(D,\,B))}
+\tau e^{-\tau(T+\text{dist}\,(\Omega,\,B))}).
\end{array}
\tag {4.8}
$$
\endproclaim

{\it\noindent Proof.}
Using one space-dimensional version of (2.14), (2.31) and (2.32),
we have
$$\begin{array}{c}
\displaystyle
I(\tau)
=-(v'(a)+c(\tau)v(a))w(a)
+O(e^{-\tau(T+\text{dist}\,(D,\,B))}
+\tau e^{-\tau(T+\text{dist}\,(\Omega,\,B))}).
\end{array}
\tag {4.9}
$$
Note that $\text{dist}\,(D,B)=a-b$ and $\text{dist}\,(\Omega,B)=-b$.

Since
$$\displaystyle
v'(x)
=-\frac{1}{2}
\left(\int_{-\infty}^xe^{-\tau(x-y)}f(y)dy
-\int_x^{\infty}e^{-\tau(y-x)}f(y)dy\right),
$$
this together with (4.5) and (4.7) gives
$$\begin{array}{c}
\displaystyle
v(a)=\frac{1}{2\tau}\int_{b-\epsilon}^be^{-\tau(a-y)}f(y)dy,\\
\\
\displaystyle
v'(a)=-\frac{1}{2}\int_{b-\epsilon}^be^{-\tau(a-y)}f(y)dy
\end{array}
$$
and thus
$$\displaystyle
v'(a)+c(\tau)v(a)
=\frac{1}{2}\left(\frac{c(\tau)}{\tau}-1\right)e^{-\tau(a-b)}
\int_{b-\epsilon}^b e^{-\tau(b-y)}f(y)dy.
\tag {4.10}
$$

Next we compute $w(a)$.
Since $w\in L^2(-\infty, a)$ one can write
$$\displaystyle
w(x,\tau)
=\frac{1}{2\tau}
\int_{-\infty}^{a}e^{-\tau\vert x-y\vert}(f(y)-e^{-\tau T}F(y,\tau))dy
+Ae^{\tau(x-a)},
\tag {4.11}
$$
where $A$ is an undetermined constant.
From this we have
$$\begin{array}{c}
\displaystyle
w'(x,\tau)
=-\frac{1}{2}
\left(\int_{-\infty}^xe^{-\tau(x-y)}(f(y)-e^{-\tau T}F(y,\tau))dy
-\int_x^{a}e^{-\tau(y-x)}(f(y)-e^{-\tau T}F(y,\tau))dy\right)\\
\\
\displaystyle
+\tau Ae^{\tau(x-a)}.
\end{array}
\tag {4.12}
$$
Substituting (4.11) and (4.12) into (4.3), one gets
$$\begin{array}{c}
\displaystyle
-\frac{1}{2}
\int_{-\infty}^ae^{-\tau(a-y)}(f(y)-e^{-\tau T}F(y,\tau))dy
+\tau A\\
\\
\displaystyle
+
c(\tau)\left(\frac{1}{2\tau}
\int_{-\infty}^{a}e^{-\tau(a-y)}(f(y)-e^{-\tau T}F(y,\tau))dy
+A\right)
+e^{-\tau T}g(\tau)=0
\end{array}
$$
From this we obtain
$$\begin{array}{c}
\displaystyle
A=-\frac{1}{2(c(\tau)+\tau)}
\left(\frac{c(\tau)}{\tau}-1\right)
\int_{-\infty}^ae^{-\tau(a-y)}f(y)dy\\
\\
\displaystyle
+\frac{e^{-\tau T}}{2(c(\tau)+\tau)}
\left(\frac{c(\tau)}{\tau}-1\right)\int_{-\infty}^ae^{-\tau(a-y)}F(y,\tau)dy
-\frac{e^{-\tau T}}{c(\tau)+\tau}g(\tau)
\end{array}
$$
and (4.11) yields
$$\begin{array}{c}
\displaystyle
w(a,\tau)\\
\\
\displaystyle
=\frac{1}{c(\tau)+\tau}\int_{-\infty}^ae^{-\tau(a-y)}f(y)dy
-\frac{e^{-\tau T}}{c(\tau)+\tau}
\left(\int_{-\infty}^ae^{-\tau(a-y)}F(y,\tau)dy
+g(\tau)\right).
\end{array}
\tag {4.13}
$$
From (4.4), we see that the integral of the
second term of (4.13) has the bound $O(\tau^{1/2})$.
From this together with the trivial equality
$$\displaystyle
\int_{-\infty}^ae^{-\tau(a-y)}f(y)dy=
e^{-\tau(a-b)}\int_{-\infty}^ae^{-\tau(b-y)}f(y)dy
$$
we obtain
$$\displaystyle
w(a,\tau)
=
\frac{e^{-\tau(a-b)}}{c(\tau)+\tau}\int_{-\infty}^ae^{-\tau(b-y)}f(y)dy
+O(e^{-\tau T}\tau^{-1/2}).
\tag {4.14}
$$
Now a combination of (4.5), (4.6), (4.9), (4.10) and (4.14) gives (4.8).

\noindent
$\Box$

Here we introduce another important assumption on $f$ which means, implicitly,
that $x=b$ is not a zero point of $f(x)$
with infinite order:
$$\displaystyle
\exists\mu\in\Bbb R\,\,
\liminf_{\tau\longrightarrow\infty}\tau^{\mu}
\left\vert\int_{b-\epsilon}^be^{-\tau(b-y)}f(y)dy\right\vert>0.
\tag {4.15}
$$

Having (4.8), one can easily obtain the following theorem.

\proclaim{\noindent Theorem 4.1.}
Let $f$ satisfy (4.5) and (4.15).
Let $T>2\text{dist}\,(D,\,B)-\text{dist}\,(\Omega,\,B)$.

(i)  Assume that $\gamma\not=1$.  We have:

if $\gamma<1$, then there exists a $\tau_0>0$ such that, for all $\tau\ge\tau_0$
$$\displaystyle
I(\tau)>0;
$$

if $\gamma>1$, then there exists a $\tau_0>0$ such that, for all $\tau\ge\tau_0$
$$\displaystyle
I(\tau)<0.
$$

In both cases we have
$$\displaystyle
\lim_{\tau\longrightarrow\infty}
\frac{1}{2\tau}
\log\vert I(\tau)\vert
=-\text{dist}\,(D,\,B).
\tag {4.16}
$$

(ii) Assume that $\gamma=1$ and $\beta\not=0$.  We have:

if $\beta<0$, then there exists a $\tau_0>0$ such that, for all $\tau\ge\tau_0$
$$\displaystyle
I(\tau)>0;
$$

if $\beta>0$, then there exists a $\tau_0>0$ such that, for all $\tau\ge\tau_0$
$$\displaystyle
I(\tau)<0.
$$

In both cases we have (4.16).

\endproclaim

Similarly to (2.51), one has
$$\displaystyle
\int_Bf(w-v)dy
=I(\tau)+O(\tau^{-1}e^{-\tau T}).
\tag {4.17}
$$
This connects two types of data asymptotically.

A combination of (4.8) and (4.17) yields the
following theorem for back-scattering case.

\proclaim{\noindent Theorem 4.2.}
Let $f$ satisfy (4.5) and (4.15).
Let $T>2\text{dist}\,(D,\,B)$.

(i)  Assume that $\gamma\not=1$.  We have:

if $\gamma<1$, then there exists a $\tau_0>0$ such that, for all $\tau\ge\tau_0$
$$\displaystyle
\int_Bf(w-v)dy>0;
$$

if $\gamma>1$, then there exists a $\tau_0>0$ such that, for all $\tau\ge\tau_0$
$$\displaystyle
\int_Bf(w-v)dy<0.
$$

In both cases we have
$$\displaystyle
\lim_{\tau\longrightarrow\infty}
\frac{1}{2\tau}
\log\left\vert\int_Bf(w-v)dy\right\vert
=-\text{dist}\,(D,\,B).
\tag {4.18}
$$

(ii)  Assume that $\gamma=1$ and $\beta\not=0$.  We have:

if $\beta<0$, then there exists a $\tau_0>0$ such that, for all $\tau\ge\tau_0$
$$\displaystyle
\int_Bf(w-v)dy>0;
$$

if $\beta>0$, then there exists a $\tau_0>0$ such that, for all $\tau\ge\tau_0$
$$\displaystyle
\int_Bf(w-v)dy<0.
$$

In both cases we have (4.18).

\endproclaim

Thus it is interesting to consider what happens in the case when
$\gamma=1$ and $\beta=0$.

In this case from (4.8) we see that
$e^{2\tau\text{dist}\,(D,\,B)}I(\tau)$ is exponentially decaying
as $\tau\longrightarrow\infty$ provided
$T>2\text{dist}\,(D,\,B)-\text{dist}\,(\Omega,\,B)$. Thus it seems
that one has to study more about the remainder term hidden in
(4.8) by $O(e^{-\tau(T+\text{dist}\,(D,\,B))} +\tau
e^{-\tau(T+\text{dist}\,(\Omega,\,B))})$. However, in fact, we see
that $u(x,t)$ and thus $I(\tau)$ are independent of $a$ by the
following argument.

Given $f\in L^2(\Bbb R)$ with compact support define
$$\displaystyle
u(x,t)=\frac{1}{2}\int_{x-t}^{x+t}f(y)dy.
\tag {4.19}
$$
We see that
$$\displaystyle
\partial_tu+\partial_xu=f(x+t).
$$
Thus if $\text{supp}\,f\subset]-\infty,\,a]$, then $f(a+t)=0$ for any $t>0$
and thus $u$ satisfies $-\partial_xu(a,t)-\partial_tu(a,t)=0,\,\,t>0$.
Of course $u$ satisfies the wave equation in the whole space and the initial conditions
$u(x,0)=0$, $\partial_tu(x,0)=f(x)$.  Therefore by the uniqueness of the weak solution of (4.1) we have $u(x,t)=u_f(x,t)$ for $(x,t)\in\,\,]-\infty,\,a[\times]0,\,T[$.
However it is easy to see that the right-hand side of (4.19)
is independent of $a$ provided $\text{supp}\,f\subset]-\infty,\,0[$.

Summing up, we conclude:

$\bullet$  if $\gamma=1$ and $\beta=0$, then $D$ is {\it invisible}
by the data $u_f(0,t)$ for $0<t<T$ or $u_f(x,t)$ for $(x,t)\in B\times\,]0,\,T[$
for any $T$ and $f$ with $\text{sup}\,f\subset\,]-\infty,\,0[$.

Finally we note that one can also determine $\gamma$ and $\beta$
after having known $\text{dist}\,(D,\,B)$.
More precisely we have the following formula.

\proclaim{\noindent Theorem 4.3.}
Let $f$ satisfy (4.5) and (4.15).

(i)  If $T>2\text{dist}\,(D,\,B)-\text{dist}\,(\Omega,\,B)$,
then as $\tau\longrightarrow\infty$ we have the following complete asymptotic expansion:
$$
\displaystyle
\frac{\displaystyle
e^{2\tau\text{dist}\,(D,\,B)}I(\tau)}
{\displaystyle
\left\vert\int_{b-\epsilon}^be^{-\tau(b-y)}f(y)dy\right\vert^2}
\sim
-\frac{\gamma-1}{2(\gamma+1)}\frac{1}{\tau}
-\frac{\beta}{(\gamma+1)^2}
\sum_{n=0}^{\infty}\left(-\frac{\beta}{\gamma+1}\right)^n\frac{1}{\tau^{n+2}}.
\tag {4.20}
$$

(ii) If $T>2\text{dist}\,(D,\,B)$, then as $\tau\longrightarrow\infty$ the function
$$
\displaystyle
\tau\longmapsto \frac{\displaystyle
e^{2\tau\text{dist}\,(D,\,B)}\int_Bf(w-v)dy}
{\displaystyle
\left\vert\int_{b-\epsilon}^be^{-\tau(b-y)}f(y)dy\right\vert^2}
$$
has the same asymptotic expansion as (4.20).

\endproclaim

{\it\noindent Proof.} (4.15) means that there exist a positive
constant $C$ and $\tau_0>0$ such that for all $\tau\ge\tau_0$
$$\displaystyle
\left\vert\int_{b-\epsilon}^be^{-\tau(b-y)}f(y)dy\right\vert\ge C\tau^{-\mu}.
\tag {4.21}
$$
This together with (4.8) gives
$$\begin{array}{c}
\displaystyle
\frac{\displaystyle
e^{2\tau\text{dist}\,(D,\,B)}I(\tau)}
{\displaystyle
\left\vert\int_{b-\epsilon}^be^{-\tau(b-y)}f(y)dy\right\vert^2}
=-\frac{1}{2\tau}\frac{(\gamma-1)\tau+\beta}{(\gamma+1)\tau+\beta}\\
\\
\displaystyle
+O(\tau^{2\mu}e^{-\tau(T-\text{dist}\,(D,\,B))}
+\tau^{2\mu+1}e^{-\tau(T-2\text{dist}\,(D,\,B)+\text{dist}\,(\Omega,\,B))}).
\end{array}
$$
Thus expanding the first term of this right-hand side  as $\tau\longrightarrow\infty$,
we obtain (4.20).

(ii) is a direct consequence of (4.17), (4.21) and (i).

\noindent
$\Box$

Note that from the coefficients of the first and second terms of
the right-hand side in (4.20) one gets $\gamma$ and successively
$\beta$,

\subsection{Some open problems}

\noindent $\bullet$  To study the asymptotic behaviour of
(1.21) as $\tau\longrightarrow\infty$ as done in Theorem 4.3 in
one-space dimensional case. We think that in three dimensions some
restriction of the geometry about $D$
should be imposed and have to use (2.14) or (2.52) instead of
(2.15).  For the purpose the asymptotic profile of $w$ on
$\partial D$ or $B$ has to be clarified as done in one space
dimensional case, see (4.14)

\noindent $\bullet$  The case when $\gamma=1$ and $\beta=0$ seems {\it pathological},
however, from a mathematical point of view it would be interesting to study the
asymptotic behaviour of $I(\tau)$ as $\tau\longrightarrow\infty$
in three-dimensions since in one-space dimensional case one can not extract the
distance of $B$ to $D$.

When the data is given by the back-scattering kernel in
the Lax-Phillips scattering theory, a corresponding positive
result in three dimensions has been announced in
Georgiev-Arnaoudov \cite{GA}.

Note also that Hansen \cite{H} considered a corresponding problem
in the Lax-Phillips scattering theory for a transparent obstacle
with a smooth $\alpha$ on $\overline D$ whose trace onto $\partial
D$ is $1$, however, its normal derivative on $\partial D$ is
nonzero everywhere.  This is the case when both of (B1) and (B2)
are not satisfied. He showed that the back-scattering kernel still
catches the values of the support function of $D$ as the leading
singularity. What can one say in our problem setting?

\noindent $\bullet$  Can one say something when the data are measured on another fixed ball $B'$ with $\overline B'\cap \overline D=\emptyset$
and $\overline B'\cap\overline B=\emptyset$ ? For example, what happens on the asymptotic behaviour
of the following integral involving $f$ and another function $f'$ with $\text{supp}\,f'=\overline B'$ as $\tau\longrightarrow\infty$:
$$\displaystyle
\int_{B'}f'(w_f-v_f)dx.
$$
We think that this is a model of the case when the emitter and receivers of the signal
are placed on different positions and what we can use is just a {\it single pair} of the emitter and receiver.

\noindent $\bullet$  It would be interesting to test the performance of formulae (1.7), (1.12), (1.17) an (1.20) by
using numerically simulated data.

\noindent $\bullet$  Expand the range of applications of the
method presented in \cite{IE0} and here to other {\it prototype}
inverse obstacle scattering problems for
acoustic/electromagentic/elastic waves or their couplings, etc.
In the framework of the Lax-Phillips scattring theory there is a result by
Georgiev \cite{GV} for a {\it moving obstacle}.
How can one treate the case in our framework?

$$\quad$$

\centerline{{\bf Acknowledgment}}

The author was partially supported by Grant-in-Aid for
Scientific Research (C)(No. 21540162) of Japan  Society for
the Promotion of Science.

$$\quad$$

\vskip1cm
\noindent
e-mail address

ikehata@math.sci.gunma-u.ac.jp

\end{document}